\documentclass[a4paper, 12pt, twoside]{article}
\usepackage{amsmath}
\usepackage{amsthm}
\usepackage{amsfonts}
\usepackage{pifont}
\usepackage{bbm}
\usepackage{color}
\usepackage{mathrsfs}
\usepackage{fancyhdr}
\usepackage{graphicx}

\usepackage{psfrag}
\usepackage{time}
\usepackage{txfonts}

\usepackage{stmaryrd}
\usepackage{mathrsfs}
\usepackage{txfonts}
\usepackage{amsfonts}

\usepackage{indentfirst}

\newcommand{\xiaosan}{\fontsize{15pt}{22pt}\selectfont}
\newcommand{\sihao}{\fontsize{14pt}{21pt}\selectfont}

\newcommand{\xiaosi}{\fontsize{12pt}{18pt}\selectfont}

\usepackage{curves}

\numberwithin{equation}{section}

\newtheorem{theorem}{ {Theorem}}[section]

\newtheorem{remark} {   {Remark}}[section]
\newtheorem{corollary} {  {Corollary}}[section]
\newtheorem{lemma} {  {Lemma}}[section]
\newtheorem{definition}{  {Definition}}[section]

\begin{document}
\setlength{\parindent}{2em}
\newpage
\fontsize{12}{22}\selectfont\thispagestyle{empty}
\renewcommand{\headrulewidth}{0pt}
 \lhead{}\chead{}\rhead{} \lfoot{}\cfoot{}\rfoot{}
\noindent

\title{\bf  H\"{o}lder Continuity of the Spectral Measures for     One-Dimensional Schr\"{o}dinger Operator in   Exponential Regime}
\author{{Wencai  Liu and Xiaoping Yuan*}\\
{\em\small School of Mathematical Sciences}\\
{\em\small  Fudan University}\\
{\em\small  Shanghai 200433, People's Republic of China}\\
{\small 12110180063@fudan.edu.cn}\\
{\small *Corresponding author: xpyuan@fudan.edu.cn}}
\date{}
\maketitle

\renewcommand{\baselinestretch}{1.2}
\large\normalsize
\begin{abstract}
  Avila and  Jitomirskaya  prove that the  spectral measure $\mu_{\lambda v, \alpha,x}^f$ of quasi-periodic Schr\"{o}dinger operator
  is  $1/2$-H\"{o}lder continuous    with  appropriate initial vector $f$,
if
 $\alpha $ satisfies Diophantine condition
  and
    $\lambda$ is small. In the present paper,  the conclusion is extended to that  for  all $\alpha$  with    $\beta(\alpha)<\infty$,  the spectral measure $\mu_{\lambda v, \alpha,x}^f$ is
    $1/2$-H\"{o}lder continuous  with   small $\lambda$, if $v$ is real
    analytic in   a neighbor  of $\{|\Im x|\leq C\beta\}$, where $C$ is a large absolute constant. In particular,
    the spectral measure $\mu_{\lambda , \alpha,x}^f$ of almost Mathieu operator is    $1/2$-H\"{o}lder continuous if $|\lambda|<e^{-C\beta}$
     with $C$ a large absolute constant.
\end{abstract}

\setcounter{page}{1} \pagenumbering{arabic}\topskip -0.82in
\fancyhead[LE]{\footnotesize  Introduction}
\section{\xiaosan \textbf{Introduction and the Main results}}
In the present paper,  we study  the quasi-periodic Schr\"{o}dinger operator   $H=H_{\lambda v,\alpha,x}$ on  $   \ell ^2(\mathbb{Z})$:

\begin{equation}\label{G11}
    (H_{\lambda v,\alpha,x}u)_{n}=u_{n+1}+u _{n-1}+ \lambda v(x+n\alpha)u_{n},
\end{equation}
where $v: \mathbb{T}=\mathbb{R} / \mathbb{Z}\rightarrow \mathbb{R} $ is the  potential,  $\lambda$ is the coupling, $\alpha  $ is the frequency, and $x $ is the phase. In particular,  the almost Mathieu operator (AMO) is given by (\ref{G11}) with $v(x)=2\cos(2\pi x)$,
  denoted by $ H_{\lambda,\alpha,x}$.
\par
 Below,   we always assume  $\alpha\in \mathbb{R}\backslash \mathbb{Q}$, and the  potential $v$ is  real analytic in a strip
  of the real axis.
  \par
The quasi-periodic Schr\"{o}dinger operator is not only related to
  some fundamental problems in physics $ {\cite{Las}}$, but also    is
fascinating  because of its remarkable richness of the related spectral theory.
In Barry Simon's list of Schr\"{o}dinger operator problems
for the twenty-first century $\cite{Sim} $, there are three problems   about AMO. The   problems  of  quasi-periodic Schr\"{o}dinger operator
have attracted  many   authors, for instance,   Avila-Jitomirskaya \cite{AJ1},\cite{AJ2},\cite{AJ3},   Avron-Simon \cite{Avr}, Bourgain-Goldstein-Schlag \cite{BGS1},\cite{BGS2} , Goldstein-Schlag \cite{GS1},\cite{GS2},\cite{GS3} and
 Jitomirskaya-Last\cite{JL1},\cite{JL2}.
\par
For $\lambda=0$, it is easy to verify  that   Schr\"{o}dinger operator  ($ \ref{G11} $)
 has   purely absolutely continuous spectrum ($[-2,2]$) by Fourier transform. We expect   the  property ( of  purely absolutely continuous spectrum)
 preserves  under sufficiently small perturbation, i.e., $\lambda$ is small. Usually there are two  smallness about $|\lambda|$. One is perturbative, meaning that the smallness   $|\lambda|$ depends not only on the potential $v$, but also on the frequency $\alpha$; the other  is  non-perturbative, meaning that the smallness
condition   only depends on the potential $v$, not on $\alpha$.
 \par
 Recall that averaging the spectral measure  $\mu_{\lambda v, \alpha,x}^{e_0}$ with respect to $x$  (see (\ref{G28})) yields
 the integrated density of states (IDS), whose    H\"{o}lder continuity  is critical to    the purely absolutely continuous spectrum.
   In the present  paper, we concern the  H\"{o}lder continuity of IDS, and generally,    of  the individual
 spectral measures $\mu_{\lambda v, \alpha,x}^{f}$.   In our another paper \cite{LIU}, we will investigate    the persistence of the   purely absolutely continuous spectrum under small perturbation by  the H\"{o}lder continuity of IDS  and some additional results in  \cite{A2},\cite{LIU2}.

  \par
  The following notions are essential in the study of equation  (\ref{G11}).
  \par
  We say $\alpha \in \mathbb{R}\backslash \mathbb{Q} $ satisfies a Diophantine condition $\text{DC}(\kappa,\tau)$ with $\kappa>0$ and $\tau>0$,
if
$$ ||k\alpha||_{\mathbb{R}/\mathbb{Z}}>\kappa |k|^{-\tau}  \text{ for   any } k\in \mathbb{Z}\backslash \{0\},$$
where $||x||_{\mathbb{R}/\mathbb{Z}}=\min_{\ell \in \mathbb{Z}}|x-\ell| $.
Let $\text{DC}=\cup_{\kappa>0,\tau>0}\text{DC}( \kappa,\tau)$. We say $\alpha$  satisfies   Diophantine condition, if $\alpha\in \text{DC}$.
 \par
Let
 \begin{equation}\label{G12}
   \beta(\alpha)=\limsup_{n\rightarrow\infty}\frac{\ln q_{n+1}}{q_n},
 \end{equation}
 where $ \frac{p_n}{q_n} $ is  the continued fraction approximants   to $\alpha$.
   One usually calls set $\{\alpha \in \mathbb{R}\backslash \mathbb{Q}|\;\beta(\alpha)>0\}$    exponential regime and set $\{\alpha\in \mathbb{R}\backslash \mathbb{Q}|\;\beta(\alpha)=0\}$   sub-exponential regime. Notice that the set DC is a real subset of the sub-exponential regime, i.e., $\text{DC} \subsetneqq\{\alpha:\beta(\alpha)=0\}$.
\par
 Here we would like to talk about some  history on   H\"{o}lder continuity of IDS, and generally,    of  the individual
 spectral measures $\mu_{\lambda v, \alpha,x}^{f}$.
 \par
 In $ \cite{E} $, Eliasson  treats  ($ \ref{G11} $)  as a dynamical systems problem--reducibility of associated cocycles. He
 shows that such cocycles are     reducible  for a.e. spectrum, and
gives out useful estimates for the non-reducible ones via a sophisticated  KAM-type methods, which breaks the limitations of the earlier KAM methods, for instance, the work of Dinaburg and Sinai \cite{Din}(they need exclude some parts of the spectrum).  As a result, Eliasson proves that   $H=H_{\lambda v,\alpha,x}$
has   purely absolutely continuous spectrum
 for $\alpha\in DC$   and $|\lambda|<\lambda_0(\alpha,v)$\footnote{$ \lambda_0(\ast)$ means $\lambda_0$   depends on $\ast$.  }.
  His student  Amor uses the sophisticated  KAM iteration  to  establish  the   1/2-H\"{o}lder continuity of   IDS  in a similar regime: $\alpha\in DC$   and $|\lambda|<\lambda_0(\alpha,v)$   \cite{Amo}.
   Amor's arguments also apply to quasi-periodic Schr\"{o}dinger operator in   multifrequency\footnote{Quasi-periodic Schr\"{o}dinger operator in multifrequency($k$ dimension)is given by ($H_{\lambda v,\alpha,x}u)_n=u_{n+1}+u _{n-1}+ \lambda v(x+n\alpha)u_n$, where $ v: \mathbb{T}^k= \mathbb{R}^k/ \mathbb{Z}^k\rightarrow \mathbb{R} $ is the potential.}.
     \par
    Both of  Eliasson and Amor's results are  perturbative (i.e., the smallness of $\lambda$ depends on $\alpha$). Such limitation are inherent to
       traditional KAM theory. The other  stronger results, i.e., non-perturbative  results, will be introduced next.
 \par
 Bourgain proves that for a.e. $\alpha $ and  $x$, $H=H_{\lambda v,\alpha,x}$ ($H_{\lambda,\alpha,x}$) has  purely absolutely continuous spectrum
 if $|\lambda|<\lambda_0(v)$ ($\lambda<1$).  Bourgain  approaches this by classical Aubry dulity and the sharp estimate of Green function
 in the regime of positive Lyapunov exponent \cite{B2}, \cite{BJ2}. By the way, in the regime of positive Lyapunov exponent,   he  sets up    the  H\"{o}lder continuity of   IDS   $N(E)$  by the H\"{o}lder continuity of
 Lyapunov exponent $L(E)$ and  Thouless formula (\cite{Avr}):
\begin{equation}\label{G13}
  L(E)=\int \ln|E-E'|dN(E').
\end{equation}
This is because, by Hilbert
transform and some theories  of singular integral operators, the H\"{o}lder continuity passes from
$L(E)$ to $N(E)$ \cite{GS1}.
Note that both $L(E)$ and $N(E)$ depend on $v$, $\lambda$ and $\alpha$, we sometimes drop the parameters dependence for simplicity.
  Earlier,  Goldstein and Schlag    $ \cite{GS1}$ have already obtained Bourgain's  results.    Concretely,  $L(E)$ and
  $N(E)$ are H\"{o}lder continuous  in the interval $[E_1,E_2]$ for strong Diophantine condition\footnote{We say $\alpha$ satisfies strong Diophantine condition
   if there exist some  $ \kappa>0$,  $ \tau>1 $such that  $$ ||k \alpha||_{\mathbb{R}/\mathbb{Z}} >\frac{\kappa}{|k|(\ln(1+|k|))^\tau}\text { for } k\in \mathbb{Z}\backslash \{0\}.$$}
   frequency $\alpha$ if $L(E) >0$ in $[E_1,E_2]$.
     They all  approach  their results
by  the  avalanche principle  and sharp large deviation theorems $  \cite{B2},\cite{GS1} $.
   Notice that  $L(E) >0$ when $|\lambda| $ is large in   non-perturbative regime
    by the  subharmonic  methods  (p.17, \cite{B2}).
     \par
    Here we would  point out some other meaningful  results.
    Suppose  $v$ is a
  trigonometric polynomial of degree $k_0$, and
assume positive Lyapunov exponents and Diophantine $\alpha$.
Goldstein-Schlag     \cite{GS2}    shows that $N(E)$ is
$( \frac{1}{2k_0}-\epsilon)$-H\"older continuous for any $\epsilon>0$. As for AMO,  combining with Aubry duality,
Goldstein and Schlag's arguments suggest    the IDS is  $(1/2-\epsilon)$-H\"older continuous    for all $\lambda\neq \pm 1$ and $\alpha\in DC$.
Their
  approach is via finite volume bounds, i.e.,  investigating  the eigenvalue  problem $H\phi=E\phi$ on a finite interval $[1,N]$ with Dirichlet
boundary conditions. The tools in \cite{GS2} have  been already turned out  to be an effective way in dealing with  the  quasi-periodic Schr\"{o}dinger operators, see   \cite{CH} and \cite{GS3} for example.
 Before \cite{GS2},   Bourgain   has  set up
$(1/2 -\epsilon)$-H\"older
continuity   for  AMO  with   $\alpha\in DC$ and  large (small)
$ |\lambda|$  perturbatively   $  \cite{B3} $.
  \par
 Avila and Jitomirskaya  address this issue by firstly  developing the quantitative version of  Aubry duality (\S 1.1).
They establish the $ 1/2 $-H\"older continuity of IDS     if  $\lambda \neq \pm1$    and $\alpha\in DC$ for AMO  \cite{AJ2}.
\par
Furthermore,  Avila proves that  $N(E)$ is $ 1/2 $-H\"older continuous with small $|\lambda|$ non-perturbatively in
sub-exponential regime (i.e., $\alpha$ satisfies
$\beta(\alpha)=0$), and   for AMO,  $N(E)$ is $ 1/2 $-H\"older continuous for all
  $ \lambda\neq \pm1$  in
sub-exponential regime.
   Note that Avila and Jitomirskaya use the quantitative version of  Aubry duality to obtain many   other  results
   of spectral  theory,
    for example, solving the sixth problem in \cite{Sim} entirely and the dry version of    Ten Martini Problem partly.
   We refer the reader to
   \cite{A2} and \cite{AJ2} for more discussion.
\par
Avila  and Jitomirskaya's  analysis also allows to investigate  a more delicate question: H\"{o}lder continuity of the individual
spectral measures. This is quite different from previous work. They  show that for
all $x$ and  vectors   $f \in\ell^1\cap \ell^2$,   the spectral measures $\mu_{\lambda v,\alpha,x}^f$ is 1/2-H\"{o}lder continuous
 uniformly  in  $x$, if   $\lambda$ is small  non-perturbatively
 and $\alpha\in DC$   \cite{AJ3}. Avila and Jitomirskaya approach  this by   the sharp estimate for the dynamics of
 Schr\"{o}dinger cocycles in \cite{AJ2} and some additional theories  in    $\cite{JL1},\cite{JL2},\cite{YB}$.
\par
In the present paper, we  extend  the quantitative version of  Aubry duality  to  all $\alpha$  with $\beta(\alpha)<\infty$.
 Together  with  Avila  and Jitomirskaya's arguments  in  \cite{AJ3},
  we obtain the following results.
 \begin{theorem}\label{Th11}
 For irrational number $ \alpha$ such that $ \beta(\alpha)<\infty$, if $v$ is   real
    analytic in   a neighbor  of $\{|\Im x|\leq C\beta\}$, where $C$ is a large absolute constant, then
   there exists some $\lambda_0=\lambda_0(v,\beta)>0$   such that
   $\mu_{\lambda v,\alpha,x}^f(J)\leq C(\lambda ,v,\alpha)|J|^{1/2}||f||^2_{\ell^1}$,
    for all intervals $J$ and all $x$ if $|\lambda|<\lambda_0$,
 where $\mu_{\lambda v,\alpha,x}^f$ is the associated spectral measure with $f\in \ell ^1 \cap \ell ^2$.
  In particular,  $\lambda_0=e^{-C\beta} $ for AMO.
 \end{theorem}
 \begin{remark}\label{Re12}
 If $ \beta(\alpha)=0$ and $v$ is   real
    analytic in a strip of real axis,  then  by Theorem \ref{Th11}  $\lambda_0=\lambda_0(v ) $, and  $\lambda_0=1$ for AMO. Those results are non-perturbative.
    Clearly, if $ 0<\beta (\alpha)<\infty$, the results obtained by Theorem \ref{Th11} are perturbative.
 \end{remark}
 \subsection{Quantitative  Aubry duality and Outline of  the present paper }
 In the present paper, we deal with the H\"{o}lder  continuity of individual spectral measure as the program of
 Avila  and Jitomirskaya \cite{AJ2}, \cite{AJ3}. Thus it is necessary to introduce  Avila  and Jitomirskaya's main
 contribution-quantitative  Aubry duality
 more details.
 \par
 Classical  Aubry duality (\S 2.2)  suggests  that Anderson
localization (only pure point spectrum with exponentially decaying eigenfunctions)  for the dual model $H_{\lambda v,\alpha,\theta}$
leads to reducibility for almost every
energy \cite{P2}. A more subtle duality theory is that pure point spectrum for almost every $\theta$  in the dual model allows to conclude purely absolutely continuous spectrum for almost every $x$ \cite{GJLS}. However,
localization  in general does not hold for every $\theta$ \cite{JS}.
This of course fits with the fact Schr\"{o}dinger  cocycles are not  reducible for all energies \cite{E}.
Thus the classical  Aubry duality  can not deal with all energies.
\par
It is therefore natural to introduce a weakened notion of localization
that could be expected to hold for every phase,
 and to develop some way  to link   the reducibility.
Avila and Jitomirskaya  make this idea come true. Namely, they introduce a new  concept: almost  localization
of   the dual model $ \{\hat H_{\lambda v  ,\alpha,\theta} \}_{ \theta\in \mathbb{R}}$,
which is a  kind of   weakened notion of localization,  and establish  a quantitative version of Aubry duality
that links local exponential decay of solutions to   eigenvalue problem  of
$\{H_{\lambda v ,\alpha,x}\}_{x \in  \mathbb{R}}$ (Lemma \ref{Le45}). See  \cite{A2} and \cite{AJ2} for more details.
\par
By some  sharp estimates  for the dynamics   of Schr\"{o}dinger  cocycles via the  quantitative version of Aubry duality, Avila and Jitomirskaya obtain  some results of  the  H\"{o}lder continuity  of IDS non-perturbatively. Together with
  the  dynamical reformulation
of
weyl-function and power-law subordinacy  techniques  in   \cite{JL1},\cite{JL2}, \cite{YB}, they set up the $\frac{1}{2}$-H\"{o}lder continuity of   individual
spectral measures, which we have said before.
\par
  Avila and Jitomirskaya's discussion is concentrated on sub-exponential regime.
   In \cite{LIU2}, we have   extended  the  quantitative version of Aubry duality  to  exponential regime    for AMO.
  In the present paper, we success to  generalize  the results of \cite{A2},\cite{AJ2} and  \cite{LIU2}, and set up  the  quantitative version of Aubry duality  for general  potential $v$  in  exponential regime.
  \par
  In order to get sharp estimate  for the dynamics   of Schr\"{o}dinger  cocycles via the  quantitative version of Aubry duality,
  the  priori estimate of the transfer matrix $A_n(x)$ is necessary,
  where  $A_n(x)$   is given by  ($\ref{G22}$) with $A=S_{\lambda v,E}$  and $E\in \Sigma_{\lambda v,\alpha}$(since the  spectrum of $H_{\lambda v,\alpha,x}$  is independent of $x$, we denote by $  \Sigma_{\lambda v,\alpha}$). In the present paper,   we obtain
    \begin{equation}\label{G14}
     ||A_n(x)||= e^{o(n)}
 \end{equation}
 through strip   $|\Im x|<\eta$ ($\eta$ will be specified later),
 which the proofs of     Avila\cite{A2} and  Avila-Jitomirskaya\cite{AJ2}  do
not apply and Avila  actually make the following footnote in  \cite{A2}:
\par
In the case of the almost Mathieu operator it is possible to show that we can take $\eta=\frac{-\ln |\lambda|}{2\pi}$ in (\ref{G14}).
For the generalization (i.e., general  potential $v$), it is possible to show that it is enough to choose $\eta$  in (\ref{G14}) such that
$v$ is holomorphic in a neighborhood of $\{|\Im x|\leq \eta\}$ and $\eta\leq \frac{1}{2\pi} \epsilon_1$, where
$\epsilon_1$ is the one in the strong localization estimate.
\par
We have confirmed (\ref{G14}) for the case of AMO \cite{LIU2}.
 In  \S 4, we will verify the claims for general $v$  by a new method.
 \par
The present paper is organized as follows:
\par
 In \S2, some basic notion will be introduced.
 In \S3, we obtain the strong localization estimate   of the Aubry dual model $\hat H_{\lambda v ,\alpha,\theta} $
 for all $\alpha$ with $\beta(\alpha) <\infty$.
In \S4,
we set up   the  priori estimate of the transfer matrix $A_n(x)$   in a given  strip.
 In \S5,
 we obtain a good  estimate  for the dynamics   of Schr\"{o}dinger  cocycles via the  quantitative  Aubry duality.
 In \S6, combining with  Avila-Jitomirskaya's  analysis in \cite{AJ3},  we prove Theorem \ref{Th11}.
\section{\xiaosan \textbf{ Preliminaries}}
\subsection{\xiaosan \textbf{  Cocycles, Lyapunov exponent, Reducibility}}
  Denote by $ \text{SL}(2,\mathbb{C})$ the all complex  $2\times 2 $-matrixes  with determinant 1.
We say a function
  $f\in C^{\omega}(\mathbb{R}/\mathbb{Z},  \mathbb{C})$ if $f$ is well  defined in $ \mathbb{R}/\mathbb{Z}$, i.e., $f(x+1)=f(x)$,
  and $f$ is   analytic in
    a strip of real axis. The definitions of $ \text{SL}(2,\mathbb{R})$ and  $  C^{\omega}(\mathbb{R}/\mathbb{Z},  \mathbb{R})$
    are similar to those of $ \text{SL}(2,\mathbb{C})$ and $C^{\omega}(\mathbb{R}/\mathbb{Z},\mathbb{C})$  respectively, except that the involved matrixes are real and the functions are real   analytic.
    \par
A $C^{\omega}$-cocycle in $ \text{SL}(2,\mathbb{C})$  is a pair $(\alpha,A)\in \mathbb{R} \times  C ^{\omega} (\mathbb{R}/\mathbb{Z}, \text{SL}(2,\mathbb{C})) $, where
$A  \in    C ^{\omega} (\mathbb{R}/\mathbb{Z}, \text{SL}(2,\mathbb{C}))$ means
      $A(x)\in \text{SL}(2,\mathbb{C})$      and the elements of $A$ are in $  C^{\omega}(\mathbb{R}/\mathbb{Z},  \mathbb{C})$.
       Sometimes, we say
     $ A$ a $C^{\omega}$-cocycle for short, if there is no ambiguity. Note that all functions, cocycles in the present paper are analytic in a strip of real axis.
     Thus we often   do not mention the analyticity, for instance, we say $A$  a cocycle instead of   $C^{\omega}$-cocycle.

    The Lyapunov exponent for the cocycle $A$ is given  by
 \begin{equation}\label{G21}
    L(\alpha,A)=\lim_{n\rightarrow\infty} \frac{1}{n}\int_{\mathbb{R}/\mathbb{Z}} \ln \| A_n(x)\|dx,
 \end{equation}
where
\begin{equation}\label{G22}
     A_n(x) = A(x+(n-1)\alpha)A(x+(n-2)\alpha)\cdots A(x).
\end{equation}
Clearly, $ L(\alpha,A)\geq0$ since $\det A(x)=1$.
\par
By  the  subadditivity of $ L_n(\alpha,A) $,  where $ L_n(\alpha,A)= \int_{\mathbb{R}/\mathbb{Z}} \ln \| A_n(x)\|dx$, one has
\begin{equation}\label{G23}
    L(\alpha,A)=\inf_{n } \frac{1}{n} L_n(\alpha,A).
 \end{equation}
 \par

Given two cocycles $(\alpha,A)$  and $(\alpha,A^{\prime})$, a conjugacy between them is a   cocycle
$ B \in C ^{\omega}(\mathbb{R}/\mathbb{Z},  \text{SL}(2,\mathbb{C}))$ such that
 \begin{equation}\label{G24}
   B(x+\alpha)^{-1}A(x)B(x)=A^{\prime}.
 \end{equation}
We say that cocycle $(\alpha,A)$  is reducible if it is  conjugate to a constant cocycle.

 \subsection{\xiaosan \textbf{ Schr\"{o}dinger  cocycles and classical Aubry duality}}
 We now consider  the  quasi-periodic Schr\"{o}dinger operator  $ H_{\lambda v,\alpha,x} $,
  the spectrum  of operator
 $ H_{\lambda v,\alpha, x}  $
 does not depend on $ x$,  denoted by  $  \Sigma_{\lambda v,\alpha}$. Indeed, shift is an unitary
 operator on $\ell^2(\mathbb{Z})$, thus $  \Sigma  _{\lambda v,\alpha,x}= \Sigma  _{\lambda v,\alpha,x+\alpha}$, where $  \Sigma _{\lambda v,\alpha,x}$
 is the spectrum of  $ H_{\lambda v,\alpha,x}  $. By  the minimality of $ x\mapsto x+\alpha$ and  continuity of  spectrum
  $  \Sigma  _{\lambda v,\alpha,x}$ with respect  to  $x$, the statement follows.
  \par
 Let
 $$
 S_{\lambda v,E}=
 \left(
     \begin{array}{cc}
       E- \lambda v & -1 \\
       1 & 0 \\
     \end{array}
   \right).
   $$
We call $(\alpha,S_{\lambda v,E})$  Schr\"{o}dinger cocycle. For AMO, we call  almost Mathieu cocycle, denoted by $(\alpha,S_{\lambda ,E})$.
\par
Note that, by dropping the symbol ¡°$v$¡±from a   notation, we indicate the corresponding  notation for AMO. For instance, denote by $\Sigma_{\lambda,\alpha}$ the spectrum of $H_{\lambda,\alpha,x}$.

\par
Fix Schr\"{o}dinger operator $  H_{\lambda v,\alpha,x}   $, we define the Aubry dual model by $\hat{H}= \hat{H}_{\lambda v,\alpha,\theta}  $,
\begin{equation}\label{G25}
  (\hat{H}\hat{u}) _ n = \sum _{k\in \mathbb{Z}}\lambda \hat{v}_ k \hat{u}_{n-k}+2\cos(2\pi \theta+n\alpha)\hat{u}_n,
\end{equation}
where $\hat{v}_ k$ is the Fourier coefficients of potential $v$.
In particular, for AMO,  it is easy to check that $ \hat{H} _{\lambda,\alpha,\theta}   =  \lambda H_{\lambda^{-1},\alpha,\theta}   $. If $\alpha\in\mathbb{R} \backslash  \mathbb{Q}$,  the spectrum of $ \hat{H}_{\lambda v,\alpha,\theta} $   is also $\Sigma_{\lambda v ,\alpha}$   \cite{GJLS}. Classical Aubry duality expresses an algebraic relation between the families of
operators $ \{\hat{H} _{\lambda v,\alpha,\theta}\}_{\theta\in\mathbb{R}} $    and $ \{ {H} _{\lambda v,\alpha,x}\}_{x\in\mathbb{R}} $
by  Bloch waves, i.e., if
 $ u:\mathbb{R}/\mathbb{Z}\mapsto\mathbb{C}$ is an $L^2$ function whose Fourier coefficients $\hat{u}$ satisfy
  $ \hat{H} _{\lambda v,\alpha,\theta}\hat{u}=E\hat{u}$, then

  $$U(x)= \left(
           \begin{array}{c }
             e^{2\pi i \theta }u(x) \\
             u(x-\alpha)\\
           \end{array}
         \right)
  $$
satisfies $S_{\lambda v,E}(x)\cdot U(x)=e^{2\pi i \theta}U(x+\alpha).$
\subsection{Spectral measure and the    integrated density of states}
Let $H$ be a bounded self-adjoint operator on $\ell ^2(\mathbb{Z})$. Then
  $(H-z)^{-1}$ is analytic
in  $ \mathbb{C} \backslash \Sigma( H)$, where $\Sigma( H)$ is the spectrum of $H$,
and we have for $ f\in \ell^2$
\begin{equation*}
  \Im\langle(H-z)^{-1}f,f\rangle=\Im z \cdot ||(H-z)^{-1}f||^2,
\end{equation*}
where $\langle\cdot,\cdot\rangle$ is the usual inner product in $\ell ^2(\mathbb{Z})$.
Thus
\begin{equation*}
 \phi_f(z)=   \langle(H-z)^{-1}f,f\rangle
\end{equation*}
is an analytic function in the upper half plane with $\Im \phi_f\geq 0$ ($\phi_f$ is a so-called Herglotz function).
\par
Therefore one has a representation
\begin{equation}\label{G26}
 \phi_f(z)=   \langle(H-z)^{-1}f,f\rangle=\int_{\mathbb{R}} \frac{1}{x-z}d\mu^{f}(x),
\end{equation}
where $\mu^{f} $ is the spectral measure   associated to vector $f$.
Alternatively, for any Borel set $\Omega\subseteq \mathbb{R}$,
\begin{equation}\label{G27}
    \mu^f(\Omega)=\langle\mathbb{E}(\Omega)f,f\rangle,
\end{equation}
where $\mathbb{E}(\Omega) $ is the corresponding spectral projection of $H$.
\par
   Denote by $\mu^{f}_{\lambda v ,\alpha,x} $
the spectral measure of Schr\"{o}dinger operator  $  H_{\lambda  v,\alpha,x}   $ and vector $f$ as before.
The  integrated density of states (IDS) $N_{\lambda v  ,\alpha}$ is obtained by   averaging the spectral measure  $\mu_{\lambda v , \alpha,x}^{e_0}$ with respect to $x$, i.e.,
\begin{equation} \label{G28}
  N_{\lambda v,\alpha}(E)=\int_{\mathbb{R}/\mathbb{Z}} \mu^{e_0}_{\lambda v  ,\alpha,x}(-\infty,E]dx,
\end{equation}
where
$e_0$ is the Dirac mass at $0 \in \mathbb{Z}$.

\subsection{Continued fraction expansion}
Define as usual for $0\leq\alpha<1,$
$$ a_0=0,\alpha_0=\alpha,$$
and inductively for $k>0,$
$$a_k=\lfloor \alpha_{k-1}^{-1}\rfloor, \alpha_k=\alpha_{k-1}^{-1}-a_k,$$
where $\lfloor t \rfloor$ denotes the greatest integer less than or equal $t$.
\par
We define
$$
\begin{array}{cc}
            p_0=0, & q_0=1, \\
            p_1=1,& q_1=a_1 ,
          \end{array}
$$
and inductively,
\begin{eqnarray*}
  p_k &=& a_k p_{k-1}+p_{k-2}, \\
  q_k &=& a_k q_{k-1}+q_{k-2}.
\end{eqnarray*}
Recall that     $\{q_n\}_{n\in \mathbb{N}}$ is the sequence of best denominators of irrational number $\alpha$,
since it satisifies
\begin{equation}\label{G29}
\forall 1\leq k <q_{n+1}, \| k\alpha\|_{\mathbb{R}/\mathbb{Z}}\geq ||q_n\alpha||_{\mathbb{R}/\mathbb{Z}}.
\end{equation}
Moreover, we also have the following estimate,

\begin{equation}\label{G210}
      \frac{1}{2q_{n+1}}\leq\Delta_n\triangleq \|q_n\alpha\|_{\mathbb{R}/\mathbb{Z}}\leq\frac{1}{q_{n+1}}.
\end{equation}

\section{ Strong localization estimate     }
Given   $ \theta \in\mathbb{R}$ and  $\epsilon_0>0$,  we say $k$ is an $\epsilon_0$-resonance for $\theta$ if
$ \| 2\theta-k\alpha\|_{\mathbb{R}/\mathbb{Z}}\leq e^{-\epsilon_0|k|}$ and
$\| 2\theta-k\alpha\|_{\mathbb{R}/\mathbb{Z}}=\min_{|j|\leq|k|} \| 2\theta-j\alpha\|_{\mathbb{R}/\mathbb{Z}}$.
\par
Clearly, $0\in \mathbb{Z}$ is  an $\epsilon_0$-resonance.
We order the $\epsilon_0$-resonances $0=|n_0|<|n_1|\leq|n_2|\cdots$. We say  $\theta$ is $\epsilon_0$-resonant if the
set of $\epsilon_0$-resonances is infinite. If $\theta$ is non-resonant, with the set of resonances $\{n_0,n_1,\cdots, n_{j_\theta}\}$,
we set $n_{j_\theta+1}=\infty$.

\par
 Below, unless stated otherwise,  $C$ is a large absolute  constant and $c$ is a small
 absolute  constant, which may change through the arguments, even when appear in the same formula. However, their  dependence on other parameters,
 will be explicitly indicated. For instance, we denote by  $C(\alpha)$   a large constant depending on $\alpha$.
Let  $||\cdot||$  be the Euclidean norms, and denote
 $||f||_\eta=\sup _{|\Im x|<\eta} ||f(x)| |$,  $||f||_0=\sup _{x\in \mathbb{R}} ||f(x) ||$.
 \par
 \begin{definition} \label{Def31}  Given a   self-adjoint operator $H$ on $\ell^2(\mathbb{Z})$,
  we say  $\phi$ is an extended state of  $H$, if
  $  H\phi=E\phi $
      with $\phi(0)=1 $   and  $ |\phi(k)|\leq  1+|k|$, where  $E\in \Sigma(H)$.
\end{definition}
\begin{definition} \label{Def32}
We say that  $ \hat{H}_{\lambda v,\alpha,\theta} $ is almost localized
  if there exists $ C_0>1$, $ \bar{C} >0$, $ \epsilon_0>0$ and $\epsilon_1>0$  such that for  any  extended state $\hat{u}$, i.e.,
 $\hat{H}_{\lambda v ,\alpha,\theta}\hat{u}=E\hat{u}$ satisfying  $\hat{u}_0=1$ and $|\hat{u}_k|\leq 1+|k|$,  where
 $E\in \Sigma_{\lambda v,\alpha}$,  then we have $|\hat{u}_k|\leq\bar{ C}e^{-\epsilon_1|k|}$  for $C_0|n_j|<|k|<C_0^{-1}|n_{j+1}|$, where
 set $\{n_j\}$ is the $\epsilon_0$-resonances for $\theta$.
 Sometimes, we also say $\hat{H}_{\lambda v ,\alpha,\theta}$
  satisfies a strong localization estimate with parameters $C_0$, $\epsilon_0$, $\epsilon_1$ and $\bar{C}$.
\end{definition}
The next theorem is our main work in this section.
\begin{theorem}\label{Th31}
Suppose irrational number $ \alpha$    satisfies $0<\beta(\alpha)<\infty $. Let   $\epsilon_0=C_1^2\beta$ and  $ \epsilon_1=C_1^3 \beta$,
  where $C_1$   is a  large absolute constant such  that it is   much larger than any absolute constant   $C$, $c^{-1} $ emerging in the present  paper.   There exists a absolute constant $C_2$
  such that
   if $v$  is analytic  in strip $ |\Im x|<C_2\beta$, then
  there exists $\lambda_0=\lambda_0(v,\beta)>0$ such that
  $ \hat{H}_{\lambda v,\alpha,\theta} $    satisfies a strong localization estimate with   parameters  $C_0=3$, $ \epsilon_0$, $\epsilon_1$
  and $\bar{C}= {C}(\lambda,v,\alpha)$, for all $\lambda$ with  $0< |\lambda|< \lambda_0 $. In particular, $\lambda_0=e^{-C _2\beta}$  for AMO.
\end{theorem}

In  $\cite{LIU2}$, we have obtained  Theorem $\ref{Th31}$ for AMO via  estimating Green function. For general  potential $v$, we  also use the sharp
 estimate of Green function  to prove Theorem \ref{Th31} by the methods of
  Avila-Jitomirskaya  in \cite{AJ2} or Bourgain-Jitomirskaya in \cite{BJ2}. Combining with our discussion in $\cite{LIU2}$,
one can obtain Theorem $\ref{Th31}$. Next, we will give a almost entire proof.
\par
Without loss of generality,      assume $ \lambda>0$.
Let $\check{H}_{\lambda v,\alpha,\theta}\triangleq\frac{1}{\lambda}\hat{H}_{\lambda v,\alpha,\theta}$,
it suffices to  prove $\check{H}_{\lambda  v ,\alpha,\theta}$ is almost localized.
  We will sometimes drop the  $E,\lambda,\alpha,\theta$-dependence   from  the notations if there is no ambiguity.
Define $H_I=R_I\check{H}R_I$, where $R_I=$ coordinate restriction to $I=[x_1,x_2]\subset\mathbb{Z}$, and denote by
$ {G}_I =(\check{H}_I-E)^{-1}$ the associated Green  function, if   $\check{H}_I-E$  is invertible. Denote  by $ {G}_{I}(x,y)$   the matrix elements of the Green  function
${G}_I$.
\par
 Assume  $\phi$ is an extended state of $\check{H}_{\lambda v,\alpha,\theta}$. Our objective is to show that  $|\phi(k)|\leq C(\lambda,v,\alpha) e^{-\epsilon_1|k|}$
 for $3|n_j|< |k|<\frac{1}{3}|n_{j+1}|$.

 It is easy to check that (p.4, \cite{BJ2})
 \begin{equation}\label{G31}
 \phi(x)=-\sum_{y\in I,k \notin I}G_I(x,y)\hat{v}_{y-k}\phi(k),
\end{equation}
for $x\in I$.
 \par
  Set $a_k=\sum_{|j|\geq|k|,jk\geq0} |j\hat{v}_j|$.
\begin{definition}\label{Def33}
Fix $m > 0$. A point $x\in\mathbb{Z}$ will be called $(m,N)$-regular if there exists an
interval $[x_1+1 ,x_2 -1]$  with $x_2=x_1+N+1$,  containing $x$   such that
\begin{equation}\label{G32}
  \sum_{y\in I, i=1,2}| G_{I}(x,y )a_{y-x_i}|<e^{-mN}  \text { for  }i=1,2;
\end{equation}
otherwise, $y$ will be called $(m,N)$-singular.
\end{definition}
  \begin{lemma}\label{Le32}
For any $m>0$, 0 is $ (m,N)$-singular   if $N>N(m) $\footnote{$N>N(m) $ means $N$ is large enough depending on $m$.}.
 \end{lemma}
\textbf{ Proof:} Otherwise,  0 is $ (m,N)$-regular,
i.e.,
there exists an
interval $[x_1+1,x_2-1]$ with $x_2=x_1+N+1$,  containing $0$   such that
\begin{equation}\label{G33}
  \sum_{y\in I, i=1,2}| G_{I}(0,y )a_{y-x_i}|<e^{-mN}   \text { for  } i=1,2.
\end{equation}
 In (\ref{G31}), let $x=0$ and  recall that $ |\phi(k)|\leq 1+|k|$, then
 \begin{eqnarray}
  \nonumber
   | \phi(0) |&=& |\sum_{y\in I,k \notin I}G_I(0,y)\hat{v}_{y-k}\phi(k) |\\
    \nonumber
     &\leq&  \sum_{y\in I,k \notin I }|G_I(0,y)\hat{v}_{y-k}|(1+|k|)  \\
      \nonumber
     &\leq& 2N\sum_{y\in I,k \notin I }|G_I(0,y)\hat{v}_{y-k}| |y-k|  \\
     \nonumber
      &\leq&  2N\sum_{y\in I, i=1,2}| G_{I}(0,y )a_{y-x_i}| \\
      &\leq&  2Ne^{-mN}<1   \label{G34}
 \end{eqnarray}
 for $N>N(m)$, which is contradicted to the hypothesis $\phi(0)=1$.$\qed$
\par
Let us denote
$$ P_N(\theta)=\det((\check{H}_{\lambda v,\alpha,\theta}-E)|_{[0,N-1]}).$$
    \par
     Following  $ {\cite{JKS}}$,
     $P_N(\theta)$ is an even function of $ \theta+\frac{1}{2}(N-1)\alpha$  and can be written as a polynomial
      of degree $N$ in $\cos2\pi (\theta+\frac{1}{2}(N-1)\alpha )$:

     \begin{equation}\label{G35}
      P_N(\theta)=\sum _{j=0}^{N}c_j\cos^j2\pi (\theta+\frac{1}{2}(N-1)\alpha)    \triangleq  Q_N(\cos2\pi  (\theta+\frac{1}{2}(N-1)\alpha)).
     \end{equation}
 Let $A_{k,r}=\{\theta\in\mathbb{R} \;|\;Q_k(\cos2\pi   \theta   )|\leq e^{(k+1)r}\} $ with $k\in \mathbb{N}$ and $r>0$.
\begin{lemma} \label{Le33}
Suppose $\beta(\alpha),\epsilon_0$ and $\epsilon_1 $ satisfy the hypothesis  of Theorem $\ref{Th31}$. Let $C_3$  be a large absolute constant.
 There exists a absolute constant $C_2$
  such that  if $v$  is analytic  in strip $ |\Im x|<C_2\beta$, then
  there exists $\lambda_0=\lambda_0(v,\beta)>0$
such that if $y\in \mathbb{Z}$ is  $(  C_1\epsilon_1,N)$-singular, $N>N(\lambda,v,\alpha)$,
 and $x\in [y-(1-\delta) N, y-\delta N]\bigcap \mathbb{Z}$ with  $\delta \in[\frac{1}{40},\frac{1}{2})$,
we have $\theta+(x+\frac{N-1}{2})\alpha $ belongs to $ A_{N,-\ln\lambda -C_3\epsilon_0}$ for all $ \lambda\in (0,\lambda_0)$.
\end{lemma}
\textbf{ Proof:}
Otherwise, there exist   $\delta \in[\frac{1}{40},\frac{1}{2})$ and  $x\in [y-(1-\delta) N, y-\delta N]\bigcap \mathbb{Z}$ such that
$\theta+(x+\frac{N-1}{2})\alpha \notin A_{N,-\ln\lambda -C_3\epsilon_0}$.
Without loss  of generality, assume  $x=0$. Thus $\theta+\frac {N-1} {2} \alpha
\notin A_{N,-\ln\lambda-C_3 \epsilon_0 }$, that is $P_N(\theta)>\lambda^{-N} e^{-C_3\epsilon_0 N}$ by (\ref{G35}). Set
 $x_1=-1,\;x_2=N$.
It  is enough to    show
that for $y\in [x_1+1,x_2-1]=I $ with $\mbox{dist} (y,\partial [x_1,x_2])
\ge \delta N$,  one has
\begin{equation}\label{G36}
(*)=\sum_{z\in
I,\,i=1,2}|G_{I}(y,z) a_{z-x_i}|<e^{-C_1\epsilon_1N}.
\end{equation}
By Cramer's rule $G_{I}(y,z)=\frac {\mu_{y,z}} {P_N(\theta)}$, where
$\mu_{y,z}$ is the corresponding minor. Together with the estimate of $\mu_{y,z}$ in Lemma \ref{Le34} and \ref{Le35} below,
we have
\begin{equation}\label{sum1} \nonumber
    (*)  \le ( \lambda e^{C_3\epsilon_0})^{N} \sum_{n=1}^{N-1}
\sum_{  {i=1,2,} {\gamma:\,|\gamma|=n}}|\det
R_{I\backslash \gamma}( \check{H}-E)R^*_{I\backslash\gamma}|
|a_{x_i-\gamma_{|\gamma|+1}}|
\prod_{i'=1}^{n}|\hat v_{\gamma_{i'+1}-\gamma_{i'}}|\;\;\;\;\;\;\;\;\;\;\;\;\;\;
\end{equation}
\begin{equation}\label{G37}
   \;\;\;\;\;\;\;\;\;\;  \leq e^{(C_3\epsilon_0+C \|v\|_0^{1/2} \lambda^{1/2}) N} \sum_{n=1}^{N-1}
\sum_{ {i=1,2,} {\gamma:\,|\gamma|=n}}  C(v,\sigma)^{n+1}
\left (\|v\|_0+C^{-1} \lambda^{-1} \frac {(n+1)^2} {N^2} \right )^{-(n+1)}
e^{- \sigma  b(\gamma,i')},
\end{equation}
where $ \sigma>0   $ is such that
\begin{equation}\label{G38}
   |\hat v_k| \leq C (v,\sigma) e^{-2|k| \sigma}
\end{equation}
 and
$b(\gamma,i')=|\gamma_{|\gamma|+1}-x_{i}|+\sum_{i'=1}^{|\gamma|}
|\gamma_{i'+1}-\gamma_{i'}|$.
Let $G_{b,n}=\{\gamma,\, |\gamma|=n \text { and } b(\gamma,i')=b\}$,  thus

\begin{equation*}
\;\;\;\;\;\;  \;\;\;\;\;   (*) \le e^{(C_3\epsilon_0+C \|v\|_0^{1/2} \lambda^{1/2}) N} \sum_{n=1}^{N-1} \sum_b
 C(v,\sigma)^{n+1}
\left (\|v\|_0+C^{-1} \lambda^{-1} \frac {(n+1)^2} {N^2} \right )^{-(n+1)}
e^{- \sigma  b} \# G_{b,n}\;\;\;\;\;\;\;\;\;\;\;\;\;\;
\end{equation*}
\begin{equation}\label{G39}
     \leq
e^{(C_3\epsilon_0+C \|v\|_0^{1/2} \lambda^{1/2}) N} \sum_{n=1}^{N-1}
 C(v,\sigma)^{n+1}
\left (C^{-1} \lambda^{-1} \frac {(n+1)^2} {N^2} \right )^{-(n+1)}
\sum_{b,\, G_{b,n} \neq \emptyset}
e^{- \sigma  b} \binom{b}{n}.
\end{equation}

If $G_{b,n} \neq \emptyset$, then $\delta N\leq \max \{ \text{dist}(y,\partial I),n+1\}
\leq b \leq (n+1) N\leq N^2$.
By Stirling formula, setting $b=r N$, $n+1=s b$, we have $\binom {b} {n} \leq
C r N e^{\phi(s) r N}$, where $\phi(s)=-s \ln
s-(1-s)\ln(1-s)$ with $ 0<s\leq1$.  Thus we have
\begin{equation}\label{G310}
   (*)\le
e^{(C_3\epsilon_0+C \|v\|_0^{1/2} \lambda^{1/2}) N}   N^5
\sup_{  ^{0<s \leq 1}  _{\delta \leq r \leq n+1}}
\left (\frac { \lambda^{-1}} {C (v,\sigma)}
r^2 s^2 \right )^{-r s N}
e^{-\sigma r N} e^{\phi(s) r N}.
\end{equation}
To prove (\ref{G36}), it suffices to show
\begin{eqnarray}\label{G311}
\nonumber
  (**) &=& \sup_{0<s \leq 1}
C_3 \epsilon_0+C\|v\|_0^{1/2} \lambda^{1/2}+\left
(  \ln C(v,\sigma)+\ln \lambda
-2 \ln r s-\frac {\sigma} {s}+\frac {\phi(s)} {s} \right ) r s \\
    &<& -2C_1\epsilon_1,
\end{eqnarray}
for any $r \in[\delta , n+1]$.
\par
Using that $\|v\|_0 \leq \frac {  C(v,\sigma)} {\sigma}$, and that
$\phi(s)/s \leq 1 -\ln s$, one has
\begin{equation}\label{G312}
    (**) \leq C_3\epsilon_0+
\left (C r c_0^{1/2}-\frac {r} {2} \right ) \sigma+
\left (  C+\ln c_0+3 \ln \frac {\sigma} {s}-
\frac {\sigma} {2 s} \right ) r s,
\end{equation}
 where   $c_0=r^{-2} \lambda C(v,\sigma) \sigma^{-3}$.
 It is easy to verify that  $ 3 \ln \frac {\sigma} {s}-
\frac {\sigma} {2 s}\leq C $,  then
\begin{equation}\label{G313}
    (**) \leq C_3\epsilon_0+
\left (C r c_0^{1/2}-\frac {r} {2} \right ) \sigma+
\left ( C+\ln c_0
 \right ) r s.
\end{equation}
Thus to show  $ (**) \leq -2C_1\epsilon_1$, it is enough to estimate (\ref{G312}) at $r=\delta$, that is
\begin{equation}\label{G314}
    (**) \leq C_3\epsilon_0+
\left (C c_0^{1/2}-\frac {1} {2} \right ) \delta \sigma+
\left (\ln C+\ln c_0
 \right ) \delta s\leq -2C_1\epsilon_1,
\end{equation}
with $c_0=\delta^{-2} \lambda C(v,\sigma) \sigma^{-3}$.
\par
 If $v$ is analytic in $|\Im x |<C_2\beta$,
  then
  \begin{equation}\label{G315}
   |\hat v_k| \leq C(v,\sigma) e^{-2\sigma|k|  },
\end{equation}
with $\sigma=\frac{C_2\beta}{4}$.
\par
If $|\lambda|<\lambda_0(v,\beta)$ such that
 \begin{equation}\label{G316}
   C c_0^{1/2}-1/2<-1/4,   C+\ln c_0<0,
 \end{equation}
then we have

 \begin{equation}\label{G317}
   (**)<C_3 \epsilon_0-\frac{C_2}{640} \beta \leq-2C_1\epsilon_1,
 \end{equation}
 since $\delta\geq 1/40$ and $C_2$ is large enough.$\qed$
      \begin{lemma} \label{Le34}(Lemma 10, $ \cite{BJ2} $)
      \begin{equation}\label{G318}
       \mu_{y,z}=\displaystyle\sum_{\gamma}\alpha_{\gamma}\det
R_{I\backslash \gamma}( \check{H}-E)R^*_{I\backslash\gamma}
\displaystyle\prod_{i=1}^{|\gamma|}|\hat v_{\gamma_{i+1}-\gamma_i}|,
      \end{equation}
      where the sum is taken over all ordered subsets
$\gamma=(\gamma_1,\ldots,\gamma_n)$ of $I$ with $\gamma_1=y$ and
$\gamma_n=z,\;|\gamma|=n-1,$
and  $\alpha_{\gamma}\in\{-1,1\}.$
      \end{lemma}

\begin{lemma}(Lemma 5.6, $\cite{AJ2}$)$ \label{Le35}$
For any $\Lambda \subset I$ and   $N>N(\lambda,v,\alpha)$,
\begin{equation}\label{G319}
    |\det R_{I\setminus \Lambda}(\check{H}-E)R^\ast_{I \setminus\Lambda})|\leq \lambda^{-N}e^{C||v||_0^{1/2}\lambda^{1/2}N}
    (||v||_0+C^{-1}\lambda^{-1}\frac{\#\Lambda^2}{N^2})^{-\# \Lambda }.
\end{equation}
\end{lemma}

\begin{definition}
     We  say that the set $\{\theta_1, \cdots ,\theta_{k+1}\}$ is $ \xi$-uniform if
      \begin{equation}\label{G320}
        \max_{ x\in[-1,1]}\max_{i=1,\cdots,k+1}\prod_{ j=1 , j\neq i }^{k+1}\frac{|x-\cos2\pi\theta_j|}
        {|\cos2\pi\theta_i-\cos2\pi\theta_j|}<e^{k\xi}.
      \end{equation}

     \end{definition}
      \begin{lemma}\label{Le36} $(\text{Lemma 9.3 },\cite{AJ1})$
      Let $ \xi_1<\xi$. If $\theta_1,\cdots,\theta_{k+1}\in A_{k,-\ln \lambda-\xi}$,  then $\{\theta_1, \cdots ,\theta_{k+1}\}$ is not $ \xi_1$-uniform  for $ k>k(\xi,\xi_1,\lambda) $.
      \end{lemma}

       Without loss of generality, assume   $3|n_j|<k<\frac{|n_{j+1}|}{3}$.
     Select $n$ such that  $q_n\leq \frac{k}{8}<q_{n+1}$ and
       let
       $s $ be the largest positive integer satisfying
$ sq_n\leq\frac{k}{8}$. Set $I_1,I_2\subset\mathbb{Z}$ as follows
 \begin{equation}\label{G321}
    I_1=[-2sq_n+1,0] \;and \;I_2=[k-2sq_n+1,k+2sq_n] , \text{ if }n_j<0,
 \end{equation}
\begin{equation}\label{G322}
      I_1=[0, 2sq_n-1]\; and \;I_2=[k-2sq_n+1,k+2sq_n] ,  \text{ if } n_j\geq0.
\end{equation}
In either case, the total number of elements in $I_1  \cup I_2$ is $6sq_n$. Let  $\theta_ {j'} = \theta + j'\alpha$ for $ j' \in I_1  \cup I_2$.
\begin{lemma}\label{Le37}$(\text{ Lemma }3.9, \cite{LIU2})$
The set $\{\theta_ {j'}\}_{   j' \in I_1  \cup I_2}$ constructed as ($\ref{G321}$) or ($\ref{G322}$) is $C\epsilon_0$-uniform for $k>k(\alpha)$
(or equivalently $n>n(\alpha)$ ).
\end{lemma}

We can now finish the proof of\textbf{ Theorem } \ref{Th31}. By Lemma  $\ref{Le36} $ and $ \ref{Le37}$,
there exists some
$j_0\in I_1\cup I_2$ such that   $ \theta_{j_0}\notin A_{6sq_n-1,-\ln \lambda- C_3\epsilon_0}$ for some absolute constant $C_3$ ($C_3$ is larger than the absolute constant $C$ emerging in Lemma \ref{Le37}).  Notice that  $y=0$ is $(C_1\epsilon_1, N)$-singular by Lemma \ref{Le32}. If
   we let $y=0$, $N=6sq_n-1$, $\delta=\frac{99}{600}$   in Lemma $\ref{Le33}$,  then   for all $j'\in I_1$, $\theta_{j'}\in A_{6sq_n-1,-\ln \lambda-C_3\epsilon_0}$
  if   $n>n(\lambda,v,\alpha)$ (or equivalently $k>k(\lambda,v,\alpha)$ ). Let $j_0\in I_2$ be such that $\theta_{j_0}\notin A_{6sq_n-1,-\ln \lambda-C_3\epsilon_0}$.
  Again by Lemma $\ref{Le33}$, $k$ is $(C_1\epsilon_1,6sq_n-1)$-regular.
 By the proof of Lemma \ref{Le32} and noting $ sq_n\geq \frac{k}{16}$, we obtain
 \begin{equation}\label{G323}
    |\phi(k)|\leq e^{-\epsilon_1k}
 \end{equation}
for $k>k(\lambda,v,\alpha)$ and $3|n_j|<k<\frac{1}{3}|n_{j+1}|$. For $k<0$, the proof is similar.  Thus
\begin{equation*}
     |\phi(k)|\leq   e^{-\epsilon_1|k|}
\end{equation*}
if $|k|>C(\lambda,v,\alpha)$  and  $3|n_j|<|k|<\frac{1}{3}|n_{j+1}|$.
That is
 \begin{equation}\label{G324}
    |\phi(k)|\leq C(\lambda,v,\alpha) e^{-\epsilon_1|k|}
 \end{equation}
for  all $k$ with  $3|n_j|<|k|<\frac{1}{3}|n_{j+1}|$.$\qed$
\par
For frequency  $\alpha$  with $\beta(\alpha)=0$,  $\hat H_{\lambda v,\alpha,\theta}$    also satisfies  strong   localization estimate with small $\lambda$.
This has been proved by Avila and Jitomirskaya  in \cite{AJ2}.
\begin{theorem}$(\text{ Theorem }5.1, \cite{AJ2})$ \label{Th38}
Assume $v$ is real analytic in a strip of real axis and $\beta(\alpha)=0$.
There exists $\lambda_0(v)>0$ such that if
 $ 0<|\lambda|<\lambda_0$, $C_0>1$,  there exist   $\epsilon_0=\epsilon_0(v,\lambda)>0$,
 $\epsilon_1=\epsilon_1(v,\lambda,C_0)>0$
such that
$\hat H_{\lambda v,\alpha,\theta}$ satisfies strong   localization estimate
with parameters $C_0,\epsilon_0,\epsilon_1 $ and $\bar{C}= {C}(\lambda,v,\alpha,C_0)$. More precisely,  for  any  extended state $\hat{u}$ of
 $\hat{H}_{\lambda v ,\alpha,\theta}$,   we have $|\hat{u}_k|\leq C(\lambda,v,\alpha,C_0)e^{-\epsilon_1|k|}$  for all $k$  with $C_0|n_j|<|k|<C_0^{-1}|n_{j+1}|$,
 where $\{n_j\}$ is the $\epsilon_0$-resonances for $\theta$.
  In particular,
$\lambda_0=1$ for AMO.
\end{theorem}

After carefully checking the details of the proof of Theorem \ref{Th38},  we can obtain another version.

\begin{theorem} \label{Th39}
Assume $v$ is real analytic in a strip of real axis and $\beta(\alpha)=0$.
There exists $\lambda_0(v)>0$ such that if
 $ 0<|\lambda|<\lambda_0$,  there exist   $\epsilon_0=\epsilon_0(v,\lambda)>0$,
 $\epsilon_1=C_1\epsilon_0$, where $C_1$ is a large absolute constant,
such that
$\hat H_{\lambda v,\alpha,\theta}$ satisfies strong   localization estimate
with parameters $C_0=3,\epsilon_0,\epsilon_1 $ and $\bar{C}= {C}(\lambda,v,\alpha)$. More precisely,  for  any  extended state $\hat{u}$ of
 $\hat{H}_{\lambda v ,\alpha,\theta}$,   we have $|\hat{u}_k|\leq C(\lambda,v,\alpha)e^{-\epsilon_1|k|}$  for all $k$ with $3|n_j|<|k|<3^{-1}|n_{j+1}|$,
 where $\{n_j\}$ is the $\epsilon_0$-resonances for $\theta$.
  In particular,
$\lambda_0=1$ for AMO.
\end{theorem}

\section{The proof of a claim from  Avila     }
 To set up the sharp estimates  for the dynamics   of Schr\"{o}dinger  cocycles via the  quantitative version of Aubry duality,    the  priori estimate of   transfer matrix $A_n(x)$    in given strip
 is of importance, where $A_n(x)$   is given by  ($\ref{G22}$) with $A=S_{\lambda v,E}$.
\begin{theorem}\label{Th41}
Suppose $\hat{H}_{\lambda v,\alpha, \theta}$  satisfies a strong localization estimate  with parameters  $C_0>1$, $ \epsilon_0$, $\epsilon_1=2\pi \eta$ and $\bar{C}$. If
$v$ is real analytic in a neighbor  of  $\{|\Im x|\leq \eta\}$, then
  $\sup_{|\Im x|<\eta}||A_k (x)||\leq   C  (\lambda,v,\alpha,\eta,\delta)e^{\delta k}$
   for any   $\delta>0$, where
   $ A  (x)  =\left(
                                                                               \begin{array}{cc}
                                                                                 E- \lambda v(x ) & -1 \\
                                                                                 1 & 0 \\
                                                                               \end{array}
                                                                             \right)
 $ with $E\in \Sigma_{\lambda v,\alpha}$.
\end{theorem}
\begin{remark}\label{Re42}
In  footnote 5 of $\cite{A2}$, Avila  think  Theorem \ref{Th41} is right, which we have mentioned in \S 1.1. We will confirm   the statements in
this section.
\end{remark}
In this section, fix $\eta=\frac{\epsilon_1}{2\pi}$.
If we can prove that  the Lyapunov exponent is vanishing  in  the   strip $|\Im x|\leq \eta$, by   Furman's uniquely ergodic theorem, Theorem \ref{Th41}
is easy to set up (see the proof of  Theorem 4.7 in \cite{LIU2}). Thus it suffices to prove  the following lemma.
\begin{lemma}\label{Le43}
Under the hypotheses   of Theorem \ref{Th41}, let  $ \alpha\in \mathbb{R}\backslash \mathbb{Q}$ and $ -\eta\leq\epsilon \leq \eta$, then
$ L(\alpha,  \epsilon  )= 0$,
where   $ L(\alpha,  \epsilon  )=  L(\alpha,A_{ \epsilon } )$ and
$$A_{\epsilon}  =\left(
                                                                               \begin{array}{cc}
                                                                                 E- \lambda v(x+i\epsilon) & -1 \\
                                                                                 1 & 0 \\
                                                                               \end{array}
                                                                             \right)
                                                                             \text{ with } E\in \Sigma_{\lambda v,\alpha}.
 $$
\end{lemma}
Following (\ref{G23}), the   Lyapunov exponent $ L(\alpha, \epsilon ) $ is  lower semi-continuous with respect to  $ \alpha\in \mathbb{R}\backslash \mathbb{Q}$  and
  $\epsilon$,   thus it is enough to  show that, for any    $\kappa,\tau>0$, $ L(\alpha, \epsilon )= 0$ if   $\alpha\in DC(\kappa,\tau) $.
\par
In this section, $\tilde{C}  $ is a large constant and $ \tilde{c}$ is a small constant.
They   are allowed to  depend on parameters  $v$, $\lambda$, $\alpha$,   $C_0$, $\bar{C} $,
$\epsilon_0,\epsilon_1,\kappa,\tau$,    which may change through the arguments, even when appear in the same formula.
 Further dependence on other parameters,
 will be explicitly indicated. For instance, we will use $ \tilde{C}(\delta)$ for a large constant depending on $\delta$, and
 $v$, $\lambda$, $\alpha$,   $C_0$, $\bar{C} $,
$\epsilon_0,\epsilon_1,\kappa,\tau$.
 \par
For the proof of   vanishing  Lyapunov exponent, a couple of lemmata  and theorems are necessary.
\par
We will say that a trigonometrical polynomial $p : \mathbb{R}/\mathbb{Z}  \mapsto \mathbb{C} $ has essential degree at most $k$ if its
Fourier coefficients   outside an interval $I$ of length $k$ ( $k=b-a$ for $I=[a,b]$) are vanishing.
\begin{lemma}\label{Le44}($\text{Theorem } 6.1, \cite{AJ2}$  )
Let $ 1\leq r\leq\lfloor q_{n+1}/q_n\rfloor$.
If $p$ has essential degree at most $k=rq_n-1$ and $x_0\in \mathbb{R}/\mathbb{Z}$, then
\begin{equation}\label{G41}
    \| p\|_0\leq C q_{n+1}^{ Cr }\sup_{0\leq j\leq k}|p(x_0+j\alpha)|.
\end{equation}
\end{lemma}
 If $\alpha\in $ DC$(\kappa,\tau)$,  then  $    q_{n+1}\leq \frac{1}{\kappa}q_n^{\tau }$ by (\ref{G29}) and (\ref{G210}), and (\ref{G41}) becomes
\begin{equation}\label{G42}
    \| p\|_0\leq Ce^{Cr\ln q_{n+1}}\sup_{0\leq j\leq k}|p(x_0+j\alpha)|\leq C e^{o( k)}\sup_{0\leq j\leq k}|p(x_0+j\alpha)|.
\end{equation}
\par

\begin{lemma}\label{Le45}$(\text{Theorem }3.3, \cite{AJ2})$
 If $E\in \Sigma_{\lambda v,\alpha} $, then there exists $ \theta\in \mathbb{R}$ and a bounded solution  of
 $\hat{H}_{\lambda v, \alpha,\theta}\hat{u}=E\hat{u}$ with
 $\hat{u}_ 0  =1$ and $|\hat{u}_ k |\leq 1$.
\end{lemma}
Given   $ E\in \Sigma_{\lambda v,\alpha}$,   let   $\theta = \theta(E)$ and solution $\hat{u}_k$ be given by Lemma $ \ref{Le45}$, and  $\{n_j\}$
be the set of $\epsilon_0$-resonances for $\theta(E)$.

\begin{lemma}\label{Le46} ($\text{Lemma } 3.1, \cite{AJ2}$  )
    If $\alpha\in \text {DC}(\kappa,\tau)$, then $|n_{j+1}|\geq a  ||2\theta-n_j\alpha||_{\mathbb{R}/\mathbb{ Z}}^{-a }\geq a  e^{ a\epsilon_0|n_j|}$, where $a=a(\kappa,\tau)$.
  \end{lemma}
  \begin{lemma}\label{Le47}$(\text{Theorem } 2.6,\;  \cite{A2} )$\label{Th49}
  Let $U:\mathbb{R}/ \mathbb{Z}\rightarrow \mathbb{C}^2$ be analytic  in $|\Im x|<\eta$. Assume that
   $ \delta_1<||U(x)|| <\delta^{-1}_2$ for all $x$ in strip  $|\Im x|<\eta$ .
  Then there exists  $B:\mathbb{R}/ \mathbb{Z}\rightarrow \text{SL}(2,\mathbb{C} )$  analytic  in $|\Im x|<\eta$ with first column $U$ and
  $||B||_\eta\leq C\delta_1^{-2}\delta_2^{-1}(1-\ln(\delta_1\delta_2))$.
  \end{lemma}
  \begin{lemma}\label{Le48}($\text{Theorem } 6.2, \cite{AJ2}$  )
    $L(\alpha, \epsilon)=0$ for $\epsilon=0$.
  \end{lemma}

\textbf{Proof of Lemma \ref{Le43}.}
\par
Let  $$\eta_1=\sup_{\epsilon}\{ \epsilon\in[0,  \eta]\;|  L(\alpha,  \xi  )= 0 \text{ for any } |\xi|\leq \epsilon\}.$$
 By the lower semi-continuity,   $L(\alpha, \epsilon)=0$ for $|\epsilon|\leq \eta_1$. Suppose Lemma \ref{Le43} does not hold, then $\eta_1<\eta$.
  Take $3\eta_2=\eta-\eta_1$.
   Let  $n=rq_k-1<q_{k+1}$ be the maxima  with $n < \frac{1}{\tilde{C} }|n_{j+1}|$ (if $\theta$ is non-resonant,    take any  $n=rq_k>\tilde{C}(\eta_2)e^{\tilde{C}(\eta_2)|n_{j_\theta}|}$),
 and let $u^I(x)=\sum_{k\in I}\hat{u}_ k e^{2\pi i kx}$ with $I=[-[\frac{n}{2}],n-[\frac{n}{2} ]]$.
 Define     $U^I(x)= \left(
           \begin{array}{c }
             e^{2\pi i \theta } u^I(x)\\
             u^I(x-\alpha)\\
           \end{array}
         \right)
  $,
 by direct computation
  \begin{equation} \label{G43}
    AU^I(x)=e^{2\pi i\theta}U^I(x+\alpha)+e^{2\pi i\theta}\left(
                                     \begin{array}{c}
                                      g(x) \\
                                      0 \\
                                     \end{array}
                                   \right),
  \end{equation}
 and the Fourier coefficients of $g(x)$ satisfy
   \begin{equation}\label{G44}
   \hat{ g}_k=\chi_I(k)(E-2\cos2\pi(\theta+k\alpha) )\hat{u} _k -\lambda\sum\chi_I(k-j)\hat{v}_ j \hat{u}_{k-j},
  \end{equation}
  where  $ \chi_I$ is the characteristic function of $I$. Since $  \hat{H}  \hat{u} = E \hat{u}$, one also has
  \begin{equation}\label{G45}
   -\hat{ g}_k=\chi_{\mathbb{Z}\backslash I}(k)(E-2\cos2\pi(\theta+k\alpha) )\hat{u}_k-\lambda\sum\chi_{\mathbb{Z}\backslash I}(k-j)\hat{v}_ j \hat{u}_{k-j}.
  \end{equation}
  Notice that  $|\hat{u}_k|<\bar{C}e^{-2\pi \eta|k|}$ for $\frac{1}{C_0} |n_j|  < |k|<C_0|n_{j+1}|$ and $| \hat{u}_k| \leq 1$ for others.
  Thus $|\hat{u}_k|<\bar{C}e^{-2\pi \eta|k|}$ for $ \tilde{C}\ln n < |k|<\tilde{C} n$  by  Lemma \ref{Le46} and $| \hat{u}_k| \leq 1$ for all $k$.
  It is easy to check that $||g||_{\eta_1+\eta_2}\leq  \tilde{C}(\eta_2) e^{-c\eta_2 n}$ and $||U^I||_{\eta_1+\eta_2}\leq \tilde{C}(\eta_2)e^{o(n)}$,  since $v$ is analytic
in  a neighbor  of  $\{|\Im x|\leq \eta\}$.
  \par
 Fix  $ \delta=\frac{1}{ {C}_1}\eta_2$, where $C_1$ is given by Theorem \ref{Th31}. Then there exists  $\xi( \lambda,v,\alpha,\eta_1,\delta)$
 with  $0<\xi <\eta_2$
such that
\begin{equation}\label{G46}
     \sup_{|\Im  x  |<\eta_1+\xi}||A_k (x)||\leq \tilde{C}(\eta_2,\delta) e^{\delta k},
\end{equation}
since  $L(\alpha, \epsilon)=0$ for $|\epsilon|\leq \eta_1$ (Theorem 4.7, \cite{LIU2}).
 \par
Next we will prove that
  the following estimate holds,
 \begin{equation}\label{G47}
    \inf_{  |\Im  x  |<\eta_1+\xi }\|U^I(x)\|\geq \tilde{c}(\eta_2,\delta) e^{- C\delta n}.
  \end{equation}
  Otherwise, let $x_0$  with $\Im x_0=t$  and $|t|<\eta_1+\xi$ such that   $ \|U^I (x_0)\|\leq  \tilde{c}(\eta_2,\delta) e^{-C\delta n}$.
  By ($ \ref{G43}$) and   ($\ref{G46}$), $||U^I(x_0 +j\alpha)||\leq \tilde{c}(\eta_2,\delta) e^{-C\delta n} $, $0 \leq j\leq n$, since $||g||_{\eta_1+\eta_2}\leq  \tilde{C}(\eta_2) e^{-c\eta_2 n}$. This implies   $|u ^I(x_0 +j\alpha)|\leq \tilde{c}(\eta_2,\delta) e^{-C\delta n} $, $0 \leq j\leq n$.
   Thus $ \|u^I_t\|_0\leq \tilde{c}(\eta_2,\delta) e^{-C\delta n}$ by  $ (\ref{G42})$, where $u^I_t(x)=u^I(x+ti)$,  contradicting  to $\int
   _{\mathbb{R}/\mathbb{Z}}u^I_t(x)dx=1$ (since $\hat{u}_0=1$).
   \par
  Let $B(x)\in \text{SL} (2,\mathbb{C})$ be the matrix,  whose first column is  $U(x)$,   given by  Lemma \ref{Le47},    then
$||B||_{ \eta_1+\xi}\leq  \tilde{C} (\eta_2,\delta) e^{  C\delta n } $. Combining with  $(\ref{G43})$, it is easy to check that

\begin{equation}\label{G48}
     B(x+\alpha)^{-1}A(x)B(x)=
     \left(
       \begin{array}{cc}
         e^{2\pi i \theta} & 0\\
         0 &  e^{-2\pi i \theta}\\
       \end{array}
     \right)+
     \left(
       \begin{array}{cc}
         \beta_1(x) & b(x) \\
         \beta_2(x) & \beta_3{x} \\
       \end{array}
     \right),
\end{equation}
where $\|b \|_{ \eta_1+\xi}\leq  \tilde{C}(\eta_2,\delta)e^{C \delta n }$, and $\|\beta_1 \|_{ \eta_1+\xi}$, $\|\beta_2 \|_{ \eta_1+\xi}$, $\|\beta_3 \|_{ \eta_1+\xi}\leq  \tilde{C}(\eta_2,\delta)e^{-c\eta_2n}$.
Taking $ \Phi=DB(x)^{-1}$, where $D=\left(
                                      \begin{array}{cc}
                                        d & 0 \\
                                        0 & d^{-1}\\
                                      \end{array}
                                    \right)
$
 with $d =  e^{-  c\eta_2 n  }$, we get
 \begin{equation}\label{G49}
     \Phi(x+\alpha) A(x)\Phi(x)^{-1}=
     \left(
       \begin{array}{cc}
         e^{2\pi i \theta} & 0\\
         0 &  e^{-2\pi i \theta}\\
       \end{array}
     \right)+
      Q(x)
\end{equation}
 where $\|Q\|_{ \eta_1+\xi}\leq \tilde{C}(\eta_2,\delta) e^{- c\eta_2 n }$ and $\|\Phi\|_{ \eta_1+\xi}\leq \tilde{C}(\eta_2,\delta)e^{ c\eta_2 n}$. Thus
 \begin{equation}\label{G410}
    \sup_{0\leq s \leq \tilde{c}(\eta_2,\delta) e^{  c \eta_2n }}\|A_s \|_{ \eta_1+\xi} \leq \tilde{ C}(\eta_2,\delta) e^{c\eta_2n},
 \end{equation}
that is
  \begin{equation}\label{G411}
     \|A_k \|_{ \eta_1+\xi} \leq  \tilde{C}(\eta_2,\delta) k^C
 \end{equation}
 with $k=\tilde{c}(\eta_2,\delta)e^{  c\eta_2 n }$.
 It follows that  $L(\alpha, \epsilon)=0$ for any  $|\epsilon|<\eta_1+\xi$, which is contradicted to the definition of $\eta_1$.$\qed$

\section{ Sharp estimate for the dynamics of Schr\"{o}dinger cocycles  }
 In section \S4, we set up the  priori estimate of the transfer matrix $A_n(x)$    in a given strip $|\Im x|<\eta$.
 In this section, we will set up sharp estimate for the dynamics of Schr\"{o}dinger cocycles.
 \par
 We first concern the exponential regime. For $ \alpha$ with $0<\beta(\alpha)<\infty$, let  $\epsilon_0,\epsilon_1,C_0$ and $ \lambda_0(v,\beta)$ be given by Theorem \ref{Th31}.    Fix $\lambda$ with $ 0<|\lambda|<\lambda_0$. Given  $E\in \Sigma_{\lambda v,\alpha}$,
    let   $\theta = \theta(E)$ and solution $\hat{u}_k$ be given by Lemma $ \ref{Le45}$, and  $\{n_j\}$
be the set of resonances for $\theta(E)$.

 Below,
  let  $A=S_{\lambda v,E} =\left(
                                                                               \begin{array}{cc}
                                                                                 E- \lambda v    & -1 \\
                                                                                 1 & 0 \\
                                                                               \end{array}
                                                                             \right)
  $. For simplicity, set  $ h_1=C_1\beta$,     $h_2=\epsilon_0 $,  $ h=\epsilon_1$.
  \par
 Before our main work, we first give some simple facts.
\begin{lemma}$(\text{Lemma } 4.2, \cite{LIU2})$ \label{Le51}
For $|n_j|>C(\alpha)$,
\begin{equation}\label{G51}
     ||2\theta-n_j\alpha||_{\mathbb{R}/\mathbb{Z}}\geq e^{-8\beta|n_{j+1}|},
\end{equation}
in particular, $|n_{j+1}|>\frac{C_1^2}{8}|n_j|$.
\end{lemma}
\begin{lemma}$(\text{Lemma } 3.1, \cite{LIU2})$\label{Le52}
The following small divisor condition holds,
\begin{equation}\label{G52}
      ||k\alpha||_{\mathbb{R}/\mathbb{Z}}\geq c(\alpha) e^{-2\beta| k|}, \text{ for any  }k\in \mathbb{Z}\backslash\{ 0\}.
\end{equation}
\end{lemma}

\begin{lemma}\label{Le53}
For any $k$ with $|k|\leq|n_j|$ and $k\neq n_j $, the following holds,
\begin{equation}\label{G53}
      ||2\theta-k\alpha||_{\mathbb{R}/\mathbb{Z}}\geq c(\alpha) e^{-4\beta| n_j|}.
\end{equation}
\end{lemma}
\textbf{Proof:}
If $ ||2\theta-n_j\alpha||_{\mathbb{R}/\mathbb{Z}}\geq c(\alpha) e^{-4\beta| n_j|}$,
 by the definition of resonance,
 \begin{equation}\label{G54}
     ||2\theta-k\alpha||_{\mathbb{R}/\mathbb{Z}}\geq||2\theta-n_j\alpha||_{\mathbb{R}/\mathbb{Z}}\geq c(\alpha) {e^{-4\beta|n_j|}} .
 \end{equation}
 If $ ||2\theta-n_j\alpha||_{\mathbb{R}/\mathbb{Z}}\leq c(\alpha)  e^{-4\beta|n_j|}  $,
  \begin{eqnarray}
 \nonumber
    ||2\theta-k\alpha||_{\mathbb{R}/\mathbb{Z}} & \geq& ||(n_j-k)\alpha||_{\mathbb{R}/\mathbb{Z}} -  ||2\theta-n_j\alpha||_{\mathbb{R}/\mathbb{Z}}\\
   \nonumber
     &\geq& c(\alpha) {e^{-4\beta|n_j|}} -c(\alpha) {e^{-4\beta |n_j|}} \\
     &\geq& c(\alpha) {e^{-4\beta|n_j|}} , \label{G55}
 \end{eqnarray}
 where the second inequality holds by   (\ref{G52}).$\qed$
 \begin{lemma}\label{Le54}
For any $k$ with $|k|\leq  C_1|n_j|$ and $k\neq n_j $, the following holds,
\begin{equation}\label{G56}
      ||2\theta-k\alpha||_{\mathbb{R}/\mathbb{Z}}\geq   e^{- CC_1\beta| n_j|},
\end{equation}
  if $ |n_j| >C(\alpha)$.
\end{lemma}
\textbf{Proof:} By (\ref{G52})
  \begin{eqnarray}
 \nonumber
    ||2\theta-k\alpha||_{\mathbb{R}/\mathbb{Z}} & \geq& ||(n_j-k)\alpha||_{\mathbb{R}/\mathbb{Z}} -  ||2\theta-n_j\alpha||_{\mathbb{R}/\mathbb{Z}}\\
   \nonumber
     &\geq& c(\alpha) {e^{-CC_1\beta|n_j|}} -  e^{-C_1^2\beta| n_j|} \\
     &\geq&   e^{-CC_1\beta|n_j|}  , \label{G57}
 \end{eqnarray}
if $|n_j| >C(\alpha)$.$\qed$

  Fix some $n = |n_j| $ and let $N = |n_{j+1}| $ if defined, otherwise let
$N = \infty$.
Let $u(x)=u^{I_1}(x)$   and      $U^{I_1}(x)= \left(
           \begin{array}{c }
             e^{2\pi i \theta } u^{I_1}(x)\\
             u^{I_1}(x-\alpha)\\
           \end{array}
         \right)
  $ with
$I_1=[-[\frac{N}{9}], -[\frac{N}{9} ]]$  as in \S 4.
  \par
   For simplicity, denote by   $C_{\star}$($c_{\star}$) a large(small) constant depending on $\lambda,v,\alpha$.
   Clearly,
  by strong localization estimate $||U^{I_1}||_{ch_i}<C_{\star}e^{Ch_i n}$, $i=1,2$.
  \par
   Following (\ref{G43})-(\ref{G45}),
   it is easy to verify that
  \begin{equation} \label{G58}
    AU^{I_1}(x)=e^{2\pi i\theta}U^{I_1}(x+\alpha)+  \left(
                                     \begin{array}{c}
                                     g(x) \\
                                      0 \\
                                     \end{array}
                                   \right)\text{ with }  ||g||_{ch}\leq   C_{\star} e^{-   c hN} .
  \end{equation}

\begin{lemma}          \label{Le55}
For   $i=1,2$,
  \begin{equation}\label{G59}
    \inf_{|\Im x|<ch_i}\|U^{I_1}(x)\|\geq c_{\star}e^{-C \beta n}.
  \end{equation}
  \end{lemma}
 \textbf{ Proof:} Following Theorem \ref{Th41} and (\ref{G58}), we can prove the lemma. See
  Theorem  4.13 in  \cite{LIU2}   for details.$\qed$

\begin{theorem}\label{Th56}
\begin{equation}\label{G510}
  \sup_{0\leq s\leq c_{\star} e^{ c h_2n  }}||A_s||_{ch_2}\leq C_{\star} e^{C\beta n}.
\end{equation}
\end{theorem}
\textbf{Proof:} It suffices to prove that if $N<\infty$, then
\begin{equation}\label{G511}
  \sup_{0\leq s\leq c_{\star}e^{ c h_2 N }}||A_s||_{ch_2}\leq C_{\star} e^{C\beta N}.
\end{equation}
Let $B(x)\in \text{SL} (2,\mathbb{C})$ be the matrix,  whose first column is  $U^{I_1}(x)$,   given by  Theorem \ref{Le47}  with
  $\eta=ch_2$,    then
$||B||_{ch_2}\leq  C_{\star}  e^{Ch_2 n} $ by   (\ref{G59}) and a simple fact   $||U^{I_1}||_{ch_2}<C_{\star}e^{Ch_2 n}$. Combining with  $(\ref{G58})$,  one easily verifies that
\begin{equation}\label{G512}
      B (x+\alpha)^{-1}A(x) B (x)=
     \left(
       \begin{array}{cc}
         e^{2\pi i \theta} & 0\\
         0 &  e^{-2\pi i \theta}\\
       \end{array}
     \right)+
     \left(
       \begin{array}{cc}
          \beta _1(x) &  b (x) \\
          \beta _2(x) &  \beta  _3{x} \\
       \end{array}
     \right)
\end{equation}
where $\| b  \|_{ch_2}<C_{\star}e^{ Ch_2 n} $, and $\| \beta _1 \|_{ch_2}  $, $\| \beta _2 \|_{ch_2}$,
$\| \beta _3\|_{ch_2}<C_{\star}e^{-chN}$.
\par
By Lemma \ref{Le51}, $$\| b  \|_{ch_2}<C_{\star}e^{ Ch_2 n}<C_{\star}e^{ C\beta N}.$$
 Solving the following  equation (by comparing the Fourier coefficients)
\begin{equation*}
      W (x+\alpha)^{-1}  \left(\begin{array}{cc}
         e^{2\pi i \theta} & b(x)\\
         0 &  e^{-2\pi i \theta}\\
       \end{array}
     \right) W(x)=
     \left(
       \begin{array}{cc}
         e^{2\pi i \theta} & b_\ell(x)\\
         0 &  e^{-2\pi i \theta}\\
       \end{array}
     \right),
\end{equation*}
where $b_\ell=\sum _{|k|\geq N}\hat{b}_ke^{2\pi i  kx}$ and
\begin{equation}\label{G514}
   W(x)=\left(
       \begin{array}{cc}
          1&  w(x) \\
          0 & 1 \\
       \end{array}
     \right),
\end{equation} we obtain
\begin{equation}\label{G513}
    \hat{w}_k=-\hat{b} _k\frac{e^{-2\pi i \theta}}{1-e^{-2\pi i (2\theta-k\alpha)}}
\end{equation}
for $|k|<N$,    and    $\hat{w}_k=0$ for $|k|\geq N$.
By small   divisor condition (\ref{G53}) (replacing $n_j$ with $n_{j+1}$ in Lemma \ref{Le53})
\begin{equation*}
   ||2\theta-k\alpha||_{\mathbb{R}/\mathbb{Z}}>c(\alpha)e^{-C\beta N} \text { for } |k|<N,
\end{equation*}
one has  $||W||_{ch_2}<C_{\star}e^{C\beta N}$.
\par
Let $B_1(x)=BW$, noting that  $\| b  _\ell \|_{ch_2}<C_{\star}e^{ -c h_2 N}$, then  $||B_1||_{ch_2}<C_{\star}e^{C\beta N}$ and
\begin{equation}\label{G515}
      B_1 (x+\alpha)^{-1}A(x) B_1 (x)=
     \left(
       \begin{array}{cc}
         e^{2\pi i \theta} & 0\\
         0 &  e^{-2\pi i \theta}\\
       \end{array}
     \right)+
     \left(
       \begin{array}{cc}
          \beta' _1(x) &  b' (x) \\
          \beta' _2(x) &  \beta ' _3{x} \\
       \end{array}
     \right)
\end{equation}
where $\| b ' \|_{ch_2}<C_{\star}e^{ -c h_2 N}$, and $\| \beta' _1 \|_{ch_2}  $, $\| \beta' _2 \|_{ch_2}$,
$\| \beta' _3\|_{ch_2}<C_{\star}e^{-chN}$.
It follows that
\begin{equation}\label{G516}
  \sup_{0\leq s\leq  c_{\star}e^{ ch_2 N }}||A_s||_{ch_2}\leq C_{\star}e^{C\beta N}.
\end{equation}
We finish the proof.$\qed$
\begin{theorem}\label{Th57}
 There exists $ B:  \mathbb{R}/\mathbb{Z}\rightarrow SL(2,\mathbb{C})$ analytic with $|| B||_{c h_1} <C_{\star} e^{C h_1 n}$
such that
\begin{equation}\label{G517}
      B (x+\alpha)^{-1}A(x) B (x)=
     \left(
       \begin{array}{cc}
         e^{2\pi i \theta} & 0\\
         0 &  e^{-2\pi i \theta}\\
       \end{array}
     \right)+
     \left(
       \begin{array}{cc}
          \beta _1(x) &  b (x) \\
          \beta _2(x) &  \beta  _3{x} \\
       \end{array}
     \right)
\end{equation}
where $\| b  \|_{ch_1}<C_{\star}e^{-ch_2 n}$, and $\| \beta _1 \|_{ch_1}$, $\| \beta _2 \|_{ch_1}$,
$\| \beta _3\|_{ch_1}<C_{\star}e^{-chN}$.
\end{theorem}
\textbf{Proof:}
Let $B_1(x)\in \text{SL} (2,\mathbb{C})$ be the matrix,  whose first column is  $U^{I_1}(x)$,   given by  Theorem \ref{Le47}  with
  $\eta=ch_1$,    then
$||B_1||_{ch_1}\leq  C_{\star}  e^{Ch_1 n} $ and
\begin{equation}\label{G518}
      B _1(x+\alpha)^{-1}A(x) B _1(x)=
     \left(
       \begin{array}{cc}
         e^{2\pi i \theta} & 0\\
         0 &  e^{-2\pi i \theta}\\
       \end{array}
     \right)+
     \left(
       \begin{array}{cc}
          \beta' _1(x) &  b' (x) \\
          \beta' _2(x) &  \beta ' _3{x} \\
       \end{array}
     \right)
\end{equation}
where $\| b ' \|_{ch_1}<C_{\star}e^{ Ch_1 n} $, and $\| \beta' _1 \|_{ch_1}  $, $\| \beta' _2 \|_{ch_1}$,
$\| \beta '_3\|_{ch_1}<C_{\star}e^{-chN}$.
\par
Let
\begin{equation}\label{G519}
    \hat{w}_k=-\hat{b}' _k\frac{e^{-2\pi i \theta}}{1-e^{-2\pi i (2\theta-k\alpha)}}
\end{equation}
for $|k|<C_1 n $  and $ k \neq n_j $,    and    $\hat{w}_k=0$ for $|k|\geq C_1n $ or $ k =n_j $.
\par
If $n\leq C(\alpha)$,  it is easy to see  that Theorem  (\ref{Th57}) has already held  by (\ref{G518}). Thus we assume $n>C(\alpha)$ so that
the  small   divisor condition (\ref{G56}) holds,   that is
\begin{equation}\label{G520}
   ||2\theta-k\alpha||_{\mathbb{R}/\mathbb{Z}}>C(\alpha)e^{-C C_1\beta n } ,
\end{equation}
for $|k|<C_1 n $  and $ k \neq n_j $.
By (\ref{G519}) and (\ref{G520}),
we have $||W||_{ch_1}<C_{\star}e^{Ch_1 n}$,
where $w(x)=\sum_{k\in \mathbb{Z}}\hat{w}_ke^{2\pi i kx}$ and
\begin{equation}\label{G521}
   W(x)=\left(
       \begin{array}{cc}
          1&  w(x) \\
          0 & 1 \\
       \end{array}
     \right).
\end{equation}

Let $B(x)=B_1W$, then  $||B||_{ch_1}<C_{\star}e^{Ch_1 n}$ and
\begin{equation}\label{G522}
      B (x+\alpha)^{-1}A(x) B (x)=
     \left(
       \begin{array}{cc}
         e^{2\pi i \theta} & 0\\
         0 &  e^{-2\pi i \theta}\\
       \end{array}
     \right)+
     \left(
       \begin{array}{cc}
          \beta _1(x) &  b'_{n_j} (x) +b^r(x)\\
          \beta _2(x) &  \beta _3(x) \\
       \end{array}
     \right),
\end{equation}
where $\| b ^{r} \|_{ch_1}<C_{\star}e^{-c h_1C_1 n}<C_{\star}e^{-ch_2 n}$,  $ b'_{n_j}(x) = \hat{b}'_{n_j}e^{2\pi n_j i x}$, and $\| \beta _1 \|_{ch_1}$, $\| \beta _2 \|_{ch_1}$,
$\| \beta_3 \|_{ch_1}<C_{\star}e^{-ch N}$.
\par
    Thus to prove Theorem \ref{Th57}, it suffices to verify
    \begin{equation}\label{G523}
|\hat b'_{n_j}| \leq  C_{\star} e^{-c h_2n   } .
\end{equation}
      Let
\begin{equation}\label{G524}
    W'(x)=  \left(
       \begin{array}{cc}
         e^{2 \pi i \theta} & b'_{n_j}(x) \\
         0 & e^{-2 \pi i \theta}\\
       \end{array}
     \right)
\end{equation}
We can compute exactly
\begin{equation}\label{G525}
    W'_s(x)=  \left(
       \begin{array}{cc}
         e^{2 \pi i s \theta} & b_{n_j}^{'s}(x) \\
         0 & e^{-2 \pi i s\theta}\\
       \end{array}
     \right),
\end{equation}
where $|b^{'s}_{n_j}(x)|=|\hat{ b}'_{n_j}\sum_{k=0}^{s-1}e^{-2\pi ik(2\theta-n_j\alpha)}|=|\hat b'_{n_j}\frac{\sin \pi s
(2\theta-n_j\alpha)}{\sin \pi(2\theta-n_j\alpha)}| \;$ if $\sin \pi (2
\theta-n_j \alpha) \neq 0$,  and $|b^{'s}_{n_j}(x)|=s |\hat b'_{n_j}|$ otherwise.
Therefore one has
  \begin{equation}\label{G526}
\|W'_s\|_0 \geq \frac {s |\hat b'_{n_j}|} {100},
0 \leq s \leq \|2 \theta-n_j \alpha\|_{ \mathbb{R}/\mathbb{ Z}}^{-1}/10.
  \end{equation}

On the other hand,
\begin{equation}\label{G527}
\|W'_s\|_0 \leq 1+s |\hat b'_{n_j}|\leq C_{\star}(1+s) e^{Ch_1n},   s \geq 0.
\end{equation}
since $||b'||_{ch_1}<C_{\star} e^{Ch_1n}$.

Since
$A=B(x+\alpha)(W'(x)+Z(x))B(x)^{-1}$,  where
 \begin{equation}\label{G528}
    Z(x)=\left(
       \begin{array}{cc}
          \beta _1(x) & b^r(x)\\
          \beta_2(x) &  \beta _3{x} \\
       \end{array}
     \right),
 \end{equation}
 after careful computation,
 \begin{equation}\label{G529}
\|A_s\|_0\geq \|B\|_0^{-2} \left (\|W'_s\|_0-\sum_{k=1}^s
\binom {s} {k}
\|Z\|^k_0 (\max_{0 \leq j < s} \|W'_j\|_0)^{1+k} \right ),
 \end{equation}

Clearly,
$\|Z\|_0 \leq C_{\star}e^{-ch_2n}$ by the estimates of elements of $Z$, thus
\begin{equation}\label{G530}
\|A_s\|_0 \geq c_{\star} e^{-Ch_1n}(\|W'_s\|_0-C_{\star}e^{-c h_2n}),   0 \leq s \leq c_{\star} e^{ch_2 n}.
\end{equation}
Combining with (\ref{G510}),
$\|W'_s\|_0 \leq C_{\star} e^{ Ch_1n}$, $0 \leq s \leq   c_{\star} e^{c h_2n    } <\|2 \theta-n_j \alpha\|_{ \mathbb{R}/\mathbb{ Z}}^{-1}/10$.
 By (\ref{G526}), we get the estimate
\begin{equation}\label{G531}
|\hat b'_{n_j}| \leq   C_{\star}  e^{-c h_2n   } .
\end{equation}
We finish the proof.$\qed$

If  frequency  $\alpha$  satisfies $\beta(\alpha)=0$, by Theorem \ref{Th39}, there exists $\lambda_0(v)>0$
such that if  $0<|\lambda|<\lambda_0(v)$, $\hat{H}_{\lambda v,\alpha,\theta}$ satisfies a strong localization estimate with parameters $\epsilon_0(\lambda,v)$,
 $\epsilon_1(\lambda,v)$, $C_0=3,\bar{C}=C(\lambda,v,\alpha)$,   where $ \epsilon_1=C_1\epsilon_0$ with $C_1$ large enough. Let $h'=\epsilon_1$, $h'_2=\epsilon_0$ and
   $h'_1=\frac{\epsilon_0}{C_1}$.  As  the proof  of Theorem \ref{Th57},  we have the following theorem.
   In order to avoid repetition,  we omit the proof.
 \begin{theorem}  \label{Th58}
 Fix some $n=|n_j|$ and let $N=|n_{j+1}|$ if defined, otherwise let $N=\infty$. Then there exists $ B:  \mathbb{R}/\mathbb{Z}\rightarrow \text{SL}(2,\mathbb{C})$ analytic with $||B||_{ch'_1}<C_{\star}e^{Ch'_1n}$ such that
 \begin{equation}\label{G532}
      B  (x+\alpha)^{-1}A(x) B  (x)=
     \left(
       \begin{array}{cc}
         e^{2\pi i \theta} & 0\\
         0 &  e^{-2\pi i \theta}\\
       \end{array}
     \right)+
     \left(
       \begin{array}{cc}
          \beta _1(x) &  b (x) \\
          \beta _2(x) &  \beta  _3{x} \\
       \end{array}
     \right)
\end{equation}
with $\| b \|_{ch'_1}<C_{\star}e^{ -ch'_2 n} $, and $\| \beta _1 \|_{ch'_1} $, $\| \beta _2 \|_{ch'_1}$,
$\| \beta _3 \|_{ch'_1}<C_{\star}e^{-ch'N}$.
 \end{theorem}

\section{Proof of   Theorem \ref{Th11}}

 Let
$\mu_{\lambda v,\alpha,x}=\mu^{e_{-1}}_{\lambda v,\alpha,x}+\mu^{e_0}_{\lambda v,\alpha, x}$, where
$e_i$ is the Dirac mass at $i \in \mathbb{Z}$. For simplicity, sometimes we  drop some parameters dependence, for example,
   replacing $\mu_{\lambda v,\alpha,x}$ with   $\mu_{ x}$ or $\mu$.
\par
Our main theorem is:
\begin{theorem}\label{Th61}
For every $0<\varepsilon<1$ and $E\in \Sigma_{\lambda v,\alpha}$, $\mu_{ x}(E-\epsilon,E+\epsilon) \leq
C _{\star} \epsilon^{1/2}$.
\end{theorem}
 The proof of   Theorem  \ref{Th61} will be given later.    Theorem \ref{Th11} can be immediately derived from Theorem \ref{Th61}.
 \par
\textbf{ Proof of  Theorem \ref{Th11}.} Since spectral measure $\mu_{ x}$    vanishes  on $\mathbb{R}\backslash \Sigma_{\lambda v,\alpha} $,
by Theorem  \ref{Th61},
\begin{equation}\label{G61}
    \mu_x(J) \leq C _{\star} |J|^{1/2} \text{ for any interval } J\subset \mathbb{R}.
\end{equation}
Let $\sigma:\ell^2(\mathbb{Z}) \rightarrow \ell^2(\mathbb{Z})$ be the shift $f(i+1)=\sigma f(i)$,  then
$\sigma H_{\lambda v,\alpha,x} \sigma^{-1}=H_{\lambda v,\alpha,x+\alpha}$.  Thus
$\mu^{\sigma f}_{x+\alpha}=\mu^f_x$ and
$\mu^{e_k}_x=\mu^{e_0}_{x+k\alpha} \leq \mu_{x+k\alpha}$.  By (\ref {G27}),
$(\mu^f_x(J))^{1/2}$ defines a semi-norm on $\ell^2(\mathbb{Z})$.
Therefore, by the triangle inequality,
\begin{eqnarray}
\nonumber
  (\mu^f_x(J))^{1/2}  &\leq  & \sum_{k \in \mathbb{Z}} |f(k)| (\mu_{x+k\alpha}(J))^{1/2}\\
    &\leq&  C _{\star} |J|^{1/4} ||f||_{\ell ^1}.\label{G62}
\end{eqnarray}
 This implies Theorem \ref{Th11}.
 \par
 Here we list two direct  corollaries from Theorem \ref{Th11}.
  \begin{corollary}\label{Co62}
  For $\alpha\in \mathbb{R}\backslash\mathbb{Q}$ with $0<\beta(\alpha)<\infty$,
if    potential $v$ is real analytic in strip $|\Im x|<C\beta$, where $C$ is a large absolute constant, then there exists $\lambda_0=\lambda_0(v,\beta)>0$ such that
 the integrated density of states  of $H_{\lambda v,\alpha,x}$ is 1/2-H\"{o}lder for   $ |\lambda|<\lambda_0$.
 In particular, $\lambda_0=e^{-C\beta}$   for AMO.
 \end{corollary}

\begin{corollary}\label{Co63}
  If irrational number $\alpha$ satisfies $\beta(\alpha)=0$,  then for any   $v\in C^{\omega}(\mathbb{R}/\mathbb{Z},\mathbb{R})$,
   there exists $\lambda_0=\lambda_0(v )>0$ such that
   the integrated density of states  of $H_{\lambda v,\alpha,x}$ is 1/2-H\"{o}lder for   $ |\lambda|<\lambda_0$.
   In particular, $\lambda_0=1$   for AMO.
\end{corollary}
\begin{remark}
For AMO, by Aubry duality, the integrated density of states  of $H_{\lambda v,\alpha,x}$ is also 1/2-H\"{o}lder if    $ |\lambda|>e^{ C\beta}$.
 \end{remark}
\subsection{Weyl function}
We will use Weyl function to estimate spectral measure. For this reason, we give some simple facts of  Weyl function firstly.
\par
Given  $E+i\epsilon$ with $E \in \mathbb{R}$ and $\epsilon>0$,
  there exists a non-zero solution $u^{+}$
of $H_{\lambda v,\alpha,x} u^{+}=(E+i\epsilon) u^{+}$ which
is $\ell^2$ at $+\infty$. The Weyl function is given by
\begin{equation}\label{G63}
    m^{+}=-\frac {u^{+}_1} {u^{+}_0}.
\end{equation}

Let
\begin{equation}\label{G64}
    M(E+i\epsilon)=\int \frac {1} {E'-(E+i\epsilon)} d\mu(E'),
\end{equation}
where $\mu=\mu_{\lambda v,\alpha,x}=\mu^{e_{-1}}_{\lambda v,\alpha,x}+\mu^{e_0}_{\lambda v,\alpha,x}$.
 Clearly,   $ M(z)$ is a Herglotz function.
It is  immediate from the definition that
\begin{equation}\label{G65}
\Im M(E+i\epsilon)\geq \frac 1{2\epsilon} \mu(E-\epsilon,E+\epsilon).
\end{equation}

Recall the usual action of $\text{SL}(2, \mathbb{C})$,
$$ \left(
     \begin{array}{cc}
       a & b \\
       c & d\\
     \end{array}
   \right)\cdot z=\frac {az+b} {cz+d}.
$$
 We define  $z_\gamma =R_{ \gamma }z$  with $\gamma\in \mathbb{R}/\mathbb{Z}$,
where
$$
R_\gamma=
\left(
       \begin{array}{cc}
         \cos2\pi \gamma & -\sin2\pi \gamma\\
          \sin2\pi \gamma &  \cos2\pi \gamma\\
       \end{array}
     \right),
$$
and
let $\psi(z)=
\sup_\gamma |z_\gamma|.$
\begin{lemma}\label{Le64}
The following  inequality holds (p. 573, \cite{AJ3}),
\begin{equation}\label{G66}
|M(z)| \leq \psi(m^+(z))\text{ for }\Im z>0 .
\end{equation}
\end{lemma}

For  $k \in \mathbb{N}$, let

\begin{equation}\label{G67}
P_{(k)}=\sum_{j=1}^k A_{2j-1}^*(x+\alpha)
A_{2j-1}(x+\alpha).
\end{equation}

Then $P_{(k)}$ is an increasing family of
positive self-adjoint operators.  In
addition, $\|P_{(k)}\|$, $\frac {\det P_{(k)}} {\|P_{(k)}\|}$
and $\det P_{(k)}$ are all increasing positive functions of $k$. Note that $\text{tr}(A_{2j-1}^{\ast}A_{2j-1})\geq 2$, then  $\|P_{(k)}\|$ (and hence $\det P_{(k)}$) is
unbounded (since  $\text{tr} P_{(k)} \geq 2k$).

\begin{lemma}$(\text{Lemma } 4.2, \cite{AJ3})$\label{Le65}
Let $\epsilon$ be such that $\det P_{(k)}=\frac {1} {4 \epsilon^2}$,
then
\begin{equation}\label{G68}
C^{-1}<\frac {\psi(m^+(E+i\epsilon))} {2\epsilon \|P_{(k)}\|}<C.
\end{equation}
\end{lemma}
\begin{theorem}\label{Th66}
For $k \geq 1$, we have
$\|P_{(k)}\| \leq C _{\star}\|(P_{(k)})^{-1}\|^{-3}$.
\end{theorem}
The proof of Theorem \ref{Th66} will be given in the end.
\par
Set $\epsilon_k=\sqrt{\frac 1{4\det
P_{(k)}}}$, i.e., $\det P_{(k)}=\frac {1} {4 \epsilon_k^2}$.

\begin{lemma}\label{Le67}
We have $\psi(m^+(E+i\epsilon_k)) \leq C _{\star} \epsilon_k^{-1/2}$.
\end{lemma}
\textbf{Proof:}
By Theorem \ref{Th66},   $\|P_{(k)}\|=\det
P_{(k)}\|(P_{(k)})^{-1}\|<\frac
{C _{\star}}{\epsilon_k^2}\|P_{(k)}\|^{-1/3}.$
Thus $\|P_{(k)}\|\leq C _{\star}\epsilon_k^{-3/2}$ and the statement
follows from (\ref{G68}). $\qed$

\textbf{  Proof of Theorem \ref{Th61}}:
 Clearly, $\lim _{k\rightarrow \infty}\epsilon_k=0$.
Following the proof of  Theorem 4.1 in \cite{AJ3} (p. 580), $\epsilon_{k+1}>c\epsilon_k$.
Combining with  (\ref {G65}),  it is enough to show that
\begin{equation}\label{G69}
  \Im M(E+i\epsilon) \leq C _{\star} \epsilon^{-1/2}
\end{equation}
holds
 for $\epsilon=\epsilon_k.$  This    follows immediately from
 (\ref{G66}) and
 Lemma \ref{Le67}.
\subsection{ Proof of Theorem \ref{Th66}}
We give two lemmata first.
\begin{lemma}$(\text{Lemma } 4.3, \cite{AJ2})$\label{XLe68}
Let
\begin{equation*}
            T(x)=\left(
     \begin{array}{cc}
       e^{2\pi i\theta} & t(x) \\
       0 &  e^{-2\pi i\theta} \\
     \end{array}
   \right)
\end{equation*}
where $t$ has a single non-zero Fourier  coefficient, i.e., $t(x)=\hat{t}_re^{2\pi irx} $. Let $X(x)=\sum_{j=1}^{k}T_{2j-1}(x)^ \ast T_{2j-1}(x)$, then
\begin{equation}\label{G610}
   ||X||_0\approx k(1+|\hat{t}_r|^2\min\{k^2,||2\theta-r\alpha||_{\mathbb{R}/\mathbb{Z}}\}),
\end{equation}
\begin{equation}\label{G611}
   ||X^{-1}||_0^{-1}\approx k,
\end{equation}
where the notation $a\approx b $ $(a,b>0)$ denotes $C^{-1}a\leq b\leq Ca$.
\end{lemma}
\begin{lemma}$(\text{Lemma } 4.4, \cite{AJ2})$\label{XLe69}
Let  $t,T$ and $X$ be as in the Lemma \ref{XLe68}. Let    $\tilde{T}:\mathbb{R}/ \mathbb{Z}\rightarrow \text{SL}(2,\mathbb{C} )$, and  set
$\tilde{X}(x)=\sum_{j=1}^{k}\tilde{T}_{2j-1}(x)^ \ast \tilde{T}_{2j-1}(x)$.
Then
\begin{equation}\label{G612}
   ||X-\tilde{X}||_0\leq 1
\end{equation}
provided that
\begin{equation}\label{G613}
   ||T-\tilde{T}||_0\leq ck^{-2} (1+2k||t||_0)^{-2}.
\end{equation}

\end{lemma}
To prove Theorem \ref{Th66}, it  is enough  to show the following lemma holds.
\begin{lemma}\label{Le68}
For $\alpha$  with  $0<\beta(\alpha)<\infty$, then
\begin{equation}\label{G614}
     \frac{||P_k||}{||P^{-1}_k||^{-3}}\leq C_{\star},  \text{ if  }C_{\star}e ^{Ch_1 n}\leq k \leq c_{\star} e^{ch_2N}.
\end{equation}
For $\alpha$  with  $\beta(\alpha)=0$, then
\begin{equation}\label{G615}
     \frac{||P_k||}{||P^{-1}_k||^{-3}}\leq C_{\star}, \text{ if  }C_{\star}e ^{Ch'_1 n}\leq k \leq c_{\star} e^{ch'_2N}.
\end{equation}
 \end{lemma}
 \textbf{Proof:}   We only give the proof of  the case  $0 <\beta(\alpha)<\infty$,  the other case is similar. Set $\Delta >n$. Let $|r| \leq \Delta $ be such that  $||2\theta-r\alpha||=\min_{|j|\leq \Delta} ||2\theta-j\alpha||$, then $|r|\geq n$.  Following the proof of  Lemma \ref{Le53},
 \begin{eqnarray}
   ||2\theta-j\alpha||_{\mathbb{R}/\mathbb{Z}} &\geq&  c(\alpha) e^{-C\beta |r|} ,  \text{ for } |j|\leq |r|, j\neq r ,\label{G616}\\
   ||2\theta-j\alpha|| _{\mathbb{R}/\mathbb{Z}}&\geq& c(\alpha) e^{-C\beta |j|},  \text{ for } |r|< |j|  \leq \Delta.\label{G617}
 \end{eqnarray}
 Using theorem  \ref{Th57}, decompose $b = t + g + q$ so that $t$ has only the Fourier coefficient $r$, i.e., $t(x)=\hat{b}_re^{2\pi irx}$, $g$ has only the Fourier coefficients
  $j \neq r$ with $| j| \leq \Delta$ and $q$ is the rest. Then
 \begin{equation}\label{G618}
    B(x+\alpha)^{-1}A(x)B(x)=T+G+H,
 \end{equation}
 where
 $$T=\left(
       \begin{array}{cc}
         e^{2\pi i\theta} & t\\
         0 &  e^{-2\pi i\theta} \\
       \end{array}
     \right),
     G=\left(
       \begin{array}{cc}
        0 & g\\
         0 &  0\\
       \end{array}
     \right),
     H=\left(
       \begin{array}{cc}
         \beta_1 & q\\
         \beta_3 &  \beta_4 \\
       \end{array}
     \right).
 $$
 Thus
 \begin{equation}\label{G619}
    ||H||_0\leq C_{\star}e^{-ch_2n} e^{-ch_1\Delta}+C_{\star}e^{-chN}.
 \end{equation}
  Solving the following equation
\begin{equation}\label{G620}
    W(x+\alpha)^{-1} (T+G)(x)W(x)=T(x),
\end{equation}
with $W(x)=\left(
        \begin{array}{cc}
          1 & w(x) \\
          0 & 1 \\
        \end{array}
      \right)
$,
then we have
\begin{equation}\label{G621}
    \hat{w}_j=-\hat{g}_j\frac{e^{-2\pi i \theta}}{1-e^{-2\pi i (2\theta-j\alpha)}},\text{ for } j\neq r,|j|\leq \Delta,
\end{equation}
and $\hat{w}_j=0$ for others. Thus
\begin{eqnarray}
\nonumber
 ||W-id||_0 & \leq&  ||w||_0 \\
 \nonumber
   &\leq& \sum _{|j|\leq |r|}|\hat{w}_j|+ \sum _{|r|<|j|\leq\Delta}|\hat{w}_j|   \\
   & \leq&   C_{\star}e^{C\beta r-ch_2n}.\label{G622}
\end{eqnarray}

Let $\Psi = BW$,
\begin{equation}\label{G623}
    ||\Psi||_0\leq  C_{\star}e^{Ch_1 n}+ C_{\star}e^{C\beta r-ch_2n }.
\end{equation}
Let $ k_\Delta\geq 0$ be maximal such that for $ 1\leq k <k_\Delta$, if we let
\begin{equation}\label{G624}
    \tilde{T}(x)=\Psi(x+\alpha)^{-1}A(x)\Psi(x)
\end{equation}
and
\begin{equation}\label{G625}
      \tilde{X}(x)=\sum_{j=1}^{k}\tilde{T}_{2j-1}(x)^ \ast \tilde{T}_{2j-1}(x) ,X=\sum_{j=1}^{k}T_{2j-1}(x)^ \ast T_{2j-1}(x),
\end{equation}
then $$ ||X-\tilde{X}||_0\leq 1 .$$.

Notice that
\begin{equation*}
  \tilde{T}(x)-T(x)= W(x+\alpha)^{-1}H(x)W(x ),
\end{equation*}
then
\begin{equation}\label{G626}
  || \tilde{T}-T||_0\leq ||W||_0^2||H||_0.
\end{equation}
Following   Lemma  \ref{XLe69},
\begin{equation}\label{G627}
 ||W||_0^2||H||_0\geq ck^{-2}_{\Delta}(1+2k_{\Delta}|\hat{b}_r|)^{-2}\geq c_{\star}k_{\Delta}^{-4},
\end{equation}
since $|\hat{b}_r|<C_{\star} $.

Thus
\begin{eqnarray}
\nonumber
   k_\Delta &\geq& \frac{c_{\star}}{||W||_0^{\frac{1}{2}}||H||_0^{\frac{1}{4}}} \\
   \nonumber
    &\geq&  \frac{c_{\star}}{(1+ C_{\star}e^{C\beta r-ch_2|n|})( C_{\star} e^{-ch_2n} e^{-ch_1\Delta}+C_{\star}e^{-chN})} \\
    &\geq& c_{\star}\min(e^{ch_1\Delta}e^{ch_2n}, e^{-C\beta\Delta}e^{chN}) \label{G628}.
\end{eqnarray}
Notice that
\begin{equation*}
     ||P_k||\leq ||\Psi||_0^4||\tilde{X}(x+\alpha)||
\end{equation*}
and
\begin{equation*}
    ||P_k^{-1}||^{-1}\geq ||\Psi||_0^{-4}||\tilde{X}(x+\alpha)^{-1}||^{-1}.
\end{equation*}
Since $||\tilde{X}||\leq ||X||+1$ and $||\tilde{X}^{-1}||\geq ||X^{-1}||^{-1} -1$ for $1\leq k < k_\Delta$.
By Lemma \ref{XLe68},
\begin{eqnarray}
  ||P_k|| &\leq& C_{\star} k(1+|\hat{b}_r|^2k^2) (e^{Ch_1  n}+e^{C\beta r-ch_2n } ) \label{G629}\\
   ||P_k^{-1}||^{-1} &\geq& c_{\star}(e^{Ch_1 n}+e^{C\beta r-ch_2n })^{-1} k.\label{G630}
\end{eqnarray}
Thus,
\begin{equation}\label{G631}
   \frac{||P_k||}{||p^{-1}_k||^{-3}} <C_{\star}|\hat{b}_r|^2(e^{Ch_1  n}+e^{C\beta r-ch_2n })^{4}+ C_{\star}\frac{1}{k^2}(e^{Ch_1  n}+e^{C\beta r-ch_2n })\leq  C_{\star},
\end{equation}
provided that $ k\geq k_\Delta ^-$, where  $k_\Delta ^-=(e^{Ch_1  n}+e^{C\beta r-ch_2n })^{1/2}$,
since  $|\hat{b}_r|^2(e^{Ch_1  n}+e^{C\beta r-ch_2n  })^{4}<C_{\star}$ by theorem $\ref{Th57}$.
We obtain that
\begin{equation}\label{G632}
     \frac{||P_k||}{||p^{-1}_k||^{-3}}\leq C_{\star}, \text{ for }  k_\Delta^-<k<k_\Delta.
\end{equation}
In order to  prove the Lemma, we have to show that for  any $k$ with $C_{\star}e ^{Ch_1  n}\leq k \leq c_{\star}e^{ch_2N} $, there exists $\Delta> n$
such that $k_\Delta^-<k<k_\Delta$.  This is easy to satisfy by setting  $\Delta=c \frac{\ln k}{\beta}$.$\qed$

           \begin{center}
           
             \end{center}
  \end{document}